\documentclass[authoryear, 11pt]{imsart}
\RequirePackage[OT1]{fontenc}
\RequirePackage{amsthm,amsmath,natbib}
\RequirePackage{hypernat}
\startlocaldefs
\numberwithin{equation}{section}
\theoremstyle{plain}

\endlocaldefs
\usepackage{psfrag, epsfig, graphicx}
\usepackage{amsfonts,amsmath,latexsym,amssymb}
\usepackage[active]{srcltx}
\usepackage{fullpage}

\newtheorem{lemma}{Lemma}
\newtheorem{theorem}{Theorem}
\newtheorem{proposition}{Proposition}

\newtheorem{example}{Example}
\newtheorem{assumption}{Assumption}
\newtheorem{remark}{Remark}
\newtheorem{definition}{Definition}

\newtheorem{assertion}{Assertion}

\renewcommand{\kappa}{\varkappa}

\newcommand{\rd}{{\rm d}}
\newcommand{\e}{\varepsilon}

\newcommand{\cA}{{\cal A}}
\newcommand{\cB}{{\cal B}}

\newcommand{\cD}{{\cal D}}
\newcommand{\cE}{{\cal E}}

\newcommand{\cH}{{\cal H}}

\newcommand{\cK}{{\cal K}}

\newcommand{\cR}{{\cal R}}

\newcommand{\cU}{{\cal U}}
\newcommand{\cV}{{\cal V}}

\newcommand{\cX}{{\cal X}}

\newcommand{\bmu}{\boldsymbol{\mu}}

\newcommand{\blg}{\boldsymbol{\mu}}
\newcommand{\bh}{{\mathbf h}}
\newcommand{\bC}{\mathbb C}

\newcommand{\bE}{\mathbb E}
\newcommand{\bF}{\mathbb F}

\newcommand{\bL}{{\mathbb L}}

\newcommand{\bN}{{\mathbb N}}

\newcommand{\bP}{{\mathbb P}}

\newcommand{\bR}{{\mathbb R}}

\newcommand{\bT}{{\mathbb T}}

\newcommand{\mA}{\mathfrak{A}}

\newcommand{\mS}{\mathfrak{S}}

\newcommand{\ms}{\mathfrak{s}}
\newcommand{\mmp}{\mathfrak{p}}
\newcommand{\mP}{\mathfrak{P}}

\newcommand{\mh}{\mathfrak{h}}

\newcommand{\mm}{\mathfrak{m}}

\newcommand{\mH}{\mathfrak{H}}

\newcommand{\ma}{\mathfrak{a}}

\newcommand{\mV}{\mathfrak{V}}
\newcommand{\mU}{\mathfrak{U}}
\newcommand{\mD}{\mathfrak{D}}
\newcommand{\epr}{\hfill\hbox{\hskip 4pt
                \vrule width 5pt height 6pt depth 1.5pt}\vspace{0.5cm}\par}

\newcommand{\blo}{\boldsymbol{0}}

\newcommand{\blk}{\boldsymbol{k}}

\newcommand{\bmh}{\boldsymbol{\mh}}
\begin{document}
\begin{frontmatter}
\title{Some new ideas in nonparametric estimation
}
\runtitle{Ideas}

\begin{aug}
\author[t1]
{\fnms{O.V.} \snm{Lepski}
\ead[label=e1]{oleg.lepski@univ-amu.fr}}
\thankstext{t1}{This work has been carried out in the framework of the Labex Archim\`ede (ANR-11-LABX-0033) and of the A*MIDEX project (ANR-11-IDEX-0001-02), funded by the "Investissements d'Avenir" French Government program managed by the French National Research Agency (ANR).}
\runauthor{O.V. Lepski}

\affiliation{Aix--Marseille Universit\'e}

\address{Institut de Math\'ematique de Marseille\\
Aix-Marseille  Universit\'e   \\
 39, rue F. Joliot-Curie \\
13453 Marseille, France\\
\printead{e1}}
\end{aug}

\begin{abstract}
In the framework of an abstract statistical model
we discuss how to use the solution of one estimation problem ({\it Problem A})
in order to construct an estimator in another, completely different,
{\it Problem B}. As a solution of {\it Problem A} we understand a data-driven selection from a given
family of estimators  $\mathbf{A}(\mH)=\big\{\widehat{A}_\mh, \mh\in\mH\big\}$ and establishing for
the selected estimator so-called oracle inequality.
If $\hat{\mh}\in\mH$ is the selected parameter and $\mathbf{B}(\mH)=\big\{\widehat{B}_\mh, \mh\in\mH\big\}$
is an estimator's collection built in  {\it Problem B} we suggest to use the estimator $\widehat{B}_{\hat{\mh}}$.
We present very general selection rule led to selector $\hat{\mh}$ and find conditions under which the estimator
$\widehat{B}_{\hat{\mh}}$ is reasonable. Our approach is illustrated by several examples related to adaptive estimation.

\end{abstract}

\begin{keyword}[class=AMS]
\kwd[Primary ]{60E15}
\kwd[; secondary ]{62G07, 62G08}
\end{keyword}

\begin{keyword}
\kwd{adaptive estimation}
\kwd{density model}
\kwd{oracle approach}
\kwd{generalized deconvolution model}
\kwd{upper fucntion}
\end{keyword}

\maketitle

\end{frontmatter}

\section{Introduction}
\label{sec:intro}

Let $(\cX^{(n)},\mathfrak{T}^{(n)},\bP^{(n)}_f, f\in\bF)$ be the statistical experiment generated by the observation $X^{(n)}$.

Let $A:\bF\to\mS_1$ and $B:\bF\to\mS_2$  be two mappings to be estimated and $\mS_1, \mS_2$ be sets endowed with  semi-metrics  $\ell$ and $\rho$ respectively.

For any $X^{(n)}$-measurable $\mS_j$-valued map $\widetilde{Q}_j, j=1,2$, and any $q\geq 1$ introduce
$$
\cR^{q}_{A}\big[\widetilde{Q}_1, f\big]=\bE^{(n)}_f\Big[\ell\big(\widetilde{Q}_1,A(f)\big)\Big]^q;\quad \cR^{q}_{B}\big[\widetilde{Q}_2, f\big]=\bE^{(n)}_f\Big[\rho\big(\widetilde{Q}_2,B(f)\big)\Big]^q,\quad f\in\bF.
$$
Here and later $\bE^{(n)}_f$ denotes the mathematical expectation w.r.t to the $\bP^{(n)}_f$  and the number $q$ is supposed to be fixed.

The main objective of the present paper can be described as follows. Assume that the problem of estimation of $A(\cdot)$, called furthermore  {\it Problem A},
 is much easier than the estimation of $B(\cdot)$ ({\it Problem B}). We will not precise here the exact meaning of "easier", which may be the theoretical
 difficulty or computational complexity or something else. One can imagine, for instance,  that {\it Problem A} has been already solved while {\it Problem B}
 is still not. It is also important to realize that {\it Problems A} and {\it B} may have completely different natures. For example, $A(f)$ is $f$,
 where $f:\bR^d\to\bR$ is the multivariate probability density  and  $\ell=\|\cdot\|_p$ is $\bL_p$-norm on $\bR^d$, while $B(f)$ is functional, i.e. $\mS_2=\bR$.
  It can be   the value of $f$ or its derivatives at a given point, some norm of this function, the entropy functional
  or Fisher information and so on. Even if both problems have the same nature,
for instance $A(f)=B(f)=f$, the loss functions (semi-metrics $\ell$ and $\rho$) can be different.
In particular one can consider
 $\ell(\cdot)=\|\cdot\|_p$ and $\rho(\cdot)=\|\cdot\|_s$,  $p\neq s$.

The problem which we address consists in finding hypotheses under which some  elements of the solution of  {\it Problem A}  could
be used in the construction of an estimation procedure for solving {\it Problem B}. Let us discuss this approach more in detail.

 The variety of  statistical procedures developed last quarter of century deal with the data-driven selection from the particular family of estimators,
 \cite{Barron-Birge}, \cite{B-B-S}, \cite{B-M}, \cite{Bunea}, \cite{C99}, \cite{cav-tsyb}, \cite{cav-golubev}, \cite{dal08}, \cite{dev-lug97}, \cite{Gol},
 \cite{GL09}, \cite{GL11}, \cite{GL12}, \cite{lepski-kerk}, \cite{LepLev98}, \cite{nemirovski00}, \cite{rigollet-tsybakov}, \cite{Tsyb03},
 \cite{Weg} among many others. A very detailed overview on this topic can be found in the recent paper \cite{lep2015}.

Suppose that we are given by the collection of estimators $\mathbf{A}(\mH)=\big\{\widehat{A}_\mh, \mh\in\mH\big\}$,
used for the estimation of the map $A(\cdot)$,  parameterized by some parameter set $\mH$ ($\widehat{A}_\mh$ depend usually on $n$ but we will omit this
dependence in the notations).
The quality of estimation is measured by the family of  risks $\big\{\cR_{A}\big[\widehat{A}_\mh, f\big], \mh\in\mH, f\in\bF\big\}$.
Let us say that {\it Problem A} is solved if one can find $X^{(n)}$-measurable element $\hat{\mh}\in\mH$ (data-driven selector) such that
the selected estimator $\widehat{A}_{\hat{\mh}}$ satisfies so-called {\it oracle-type} inequality:
\begin{equation}
\label{eq:ell-oracle-inequality}
\cR_{A}\big[\widehat{A}_{\hat{\mh}}, f\big]\leq \inf_{\mh\in\mH}\cA_{\ell}^{(n)}\left(f,\mh\right)+cr_n,\quad\forall f\in\bF,\;\forall n\geq 1.
\end{equation}
Here $c>0$ is a numerical constant independent on $n$ and $f$ and $r_n\to 0, n\to\infty$ is given sequence.
As to the quantity
$\cA_{\ell}^{(n)}(\cdot,\cdot)$  it is explicitly  expressed and for some particular problems one can prove the inequality (\ref{eq:ell-oracle-inequality}) with
$$
\cA_{\ell}^{(n)}(f,\mh)=C\cR_{A}\big[\widehat{A}_\mh, f\big],
$$
where $C$ is as previously  a universal constant.

Let $\mathbf{B}(\mH)=\big\{\widehat{B}_\mh, \mh\in\mH\big\}$ be a given collection of statistical procedures related to the estimation of $B(\cdot)$.
It is extremely important for us that both collections $\mathbf{A}(\mH)$ and $\mathbf{B}(\mH)$ are parameterized by one and the same set $\mH$.
As it was mentioned above {\it Problem A} and {\it Problem B} may have different nature and often only the  set $\mH$ relates both problems.

\smallskip

The questions which we want to answer are: {\it under which conditions the selector $\hat{\mh}$ provides a reasonable choice from the collection $\mathbf{B}(\mH)$? Is it possible
to establish an  oracle inequality similar to (\ref{eq:ell-oracle-inequality}) for the selected estimator $\widehat{B}_{\hat{\mh}}$?}

\smallskip

We do not think that the answers on aforementioned questions can be obtained when an arbitrary selection rule led to $\hat{\mh}$ is considered. So, in the next section in the framework of an abstract statistical model  we propose a quite general selection scheme and establish for it the oracle inequality (\ref{eq:ell-oracle-inequality}). This part of the paper has an independent interest since all results will be obtained under few very general assumptions which can be verified for many statistical models and problems.
The proposed approach can be viewed as a generalization of several estimation procedures developed  by the author and his collaborators  during last twenty years, see \cite{LepLev98}, \cite{lepski-kerk}, \cite{ioud}, \cite{GL08}, \cite{GL09}, \cite{GL11}, \cite{GL12}, \cite{lepski13a}, \cite{GL13},  and \cite{lep2015}.

Coming back to the study of  the behavior of the "plug-in" estimator $\widehat{B}_{\hat{\mh}}$ we would like to emphasize that it will be done under the following assumption imposed on the statistical experiment.
We will assume that $X^{(n)}=\big(X_1^{(n)},X_2^{(n)}\big)$, where $X_1^{(n)},X_2^{(n)}$ are independent random elements. This fundamental assumption may correspond to splitting data on two independent samples or to the situation when a statistician disposes an independent copy of the considered  statistical model.
Our selection rule (led in particular to the solution of {\it Problem A}) is based on the  observation $X_1^{(n)}$ while the estimator's collection $\mathbf{B}(\mH)$ is built  from the observation $X_2^{(n)}$.

The following notations will be used in the sequel: $\bP^{(n)}_{1,f}$ and $\bP^{(n)}_{2,f}$
denote marginal laws of  $X_1^{(n)}$ and $X_2^{(n)}$ respectively and  $\bE^{(n)}_{i,f}, i=1,2$, will be used for the mathematical expectation w.r.t. $\bP^{(n)}_{i,f}$.

We finish this introduction  presenting several examples in which the discussed above strategy can be applied.
These considerations will be  continued in Sections \ref{sec:deconv} and \ref{sec:derivatives}.

\paragraph{Estimation in the deconvolution density model under $\bL_p$-loss.} Consider the following observation scheme:
\begin{equation}
\label{eq:observation-scheme}
Z_i=X_i+Y_i,\quad i=1,\ldots,n,
\end{equation}
where $X_i, i=1,\ldots,n,$ are {\it i.i.d.} $d$-dimensional  random vectors with common density $f$ to be estimated.
The noise variables  $Y_i, i=1,\ldots,n,$
are {\it i.i.d.} $d$-dimensional random vectors with known common density $g$.
The sequences  $\{X_i, i=1,\ldots,n\}$ and  $\{Y_i, i=1,\ldots,n\}$ are supposed to be mutually independent.

Let $X^{(n)}=(Z_1,\ldots, Z_n)$ and let $B(f)=f$ and $\rho(\cdot)=\ell(\cdot)=\|\cdot\|_p$ that means that we are interested in estimation of the underlying density $f$ under $\bL_p$-loss.

Let $K:\bR^d\to\bR$ be the  function belonging to $\bL_1\big(\bR^d\big)$ and  $\int_{\bR}K=1$.
Let  $\cH^d$ be the diadic grid of $(0,1]^d$ and
define for any $\vec{h}=(h_1,\ldots,h_d)\in\cH^d$
\begin{eqnarray}
\label{eq:def-kernel-K}
K_{\vec{h}}(t)=V^{-1}_{\vec{h}}K\big(t_1/h_1,\ldots,t_d/h_d\big),\; t\in\bR^d,\quad V_{\vec{h}}=\prod_{j=1}^dh_j.&&
\end{eqnarray}
For any $\vec{h}\in \cH^d$ let  $M\big(\cdot,\vec{h}\big):\bR^d\to\bR$
satisfy the operator equation
\begin{eqnarray}
\label{eq:def-kernel-M}
K_{\vec{h}}(y)=\int_{\bR^d}g(t-y)M\big(t,\vec{h}\big)\rd t,\quad  y\in\bR^d.&&
\end{eqnarray}
Introduce the following estimator's collection
\begin{equation}
\label{eq:def-M-est}
\mathbf{B}\big(\cH^d\big)=\bigg\{\widehat{B}_{\vec{h}}(\cdot)=n^{-1}\sum_{i=1}^n M\big(Z_i-\cdot,\vec{h}\big),\;\vec{h}\in\cH^d \bigg\}.
\end{equation}
In  \cite{{comte}} and in \cite{rebelles15}  the   data-driven selection rules from $\mathbf{B}\big(\cH^d\big)$, based on the methodology developed in \cite{GL12}, were proposed.
The authors established oracle inequalities of type (\ref{eq:ell-oracle-inequality}) and deduced from them several results related to adaptive estimation.

Our idea is quite different and consists in the following. Consider first the estimation of $A(f)=g\star f$, where here and later  "$\star$" denotes
the convolution operator. Note that $g\star f$ is the density of $Z_1$ and, therefore, can be easily estimated from the observation $X^{(n)}$ using standard kernel estimator. Consider the collection
\begin{equation}
\label{eq:def-K-est}
\mathbf{A}\big(\cH^d\big)=\bigg\{\widehat{A}_{\vec{h}}(\cdot)=\sum_{i=1}^nK_{\vec{h}}(Z_i-\cdot),\;\vec{h}\in\cH^d \bigg\}.
\end{equation}
The problem of bandwidth selection from $\mathbf{A}\big(\cH^d\big)$  was studied in \cite{GL11} and \cite{GL13} where several oracle inequalities were proved. Let $\vec{\mathbf{h}}_{\mathbf{n}}$ be a selected bandwidth. We propose then to use
$
\widehat{B}_{\vec{\mathbf{h}}_{\mathbf{n}}}(\cdot)
$
as the final estimator.

We remark that similar strategy in the linear  regression model (selection from the family of  spectral cut-off estimators, $d=1$)
was adopted in \cite{chern-gol}. In the sequence space gaussian white noise model the approach discussed above was applied in \cite{knapik}
in order to find the posterior contraction rate in inverse problems
in the context of the bayesian nonparametrics.


\paragraph{Estimation of derivatives in the density model under $\bL_p$-loss.} Let $X^{(n)}=(X_1,\ldots,X_n)$ be {\it i.i.d.} random variables with unknown common density $f$ (the multidimensional version of the problem will be presented in Section \ref{sec:derivatives}). Let we are interested in estimating of $f^{(m)}, m\in\bN^*$, where  $f^{(m)}$ denotes $m$-th derivative of $f$.

Thus, $B_m(f)=f^{(m)}$ and let $\rho(\cdot)=\ell(\cdot)=\|\cdot\|_p$. Set $A(f)=f$ and consider the family of kernel estimators
$$
\mathbf{A}\big(\cH^1\big)=\bigg\{\widehat{A}_{h}(\cdot)=\sum_{i=1}^nK_{h}(X_i-\cdot),\; h\in\cH^1 \bigg\}.
$$
It is worth noting that for any $m\in\bN^*$ the estimator
$
\widehat{B}_{h,m}(\cdot)=\widehat{A}^{(m)}_{h}(\cdot)
$
is usually used for estimating $B_m(f)$. Hence, one of the possibilities consists in the selection from the collection
$$
\mathbf{B}_m\big(\cH^1\big)=\big\{\widehat{B}_{h,m}(\cdot),\; h\in\cH^1 \big\}
$$
in order to construct an  estimator for $f^{(m)}$. Note, however that in this case the corresponding selector as well as the selected estimator will depend on $m$.
In particular, the selection schemes are different for different values of $m$.

Our approach to considered problem consists  in selecting from the family  $\mathbf{A}\big(\cH^1\big)$ that provides us with the selector $\mathbf{h_n}$.
Then for any $m\in\bN^*$, we suggest  to use
$
\widehat{B}_{\mathbf{h_n},m}(\cdot)
$
as an estimator for $B_m(f)$. In other words the problem we address can be formulated as follows.
Is it possible to differentiate the kernel estimator with data-driven bandwidth in order to
get the estimator for any derivative of the underlying density simultaneously?

In the framework of gaussian white noise model, this problem was studied in \cite{Efr98}.
It was shown that the answer is positive when  $p=2$ and the adaptation is considered over the collection of
Sobolev classes. More precisely it was shown that the differentiation of the Efromovich-Pinsker estimator leads to the efficient adaptive estimator
of any derivative.

In Section \ref{sec:derivatives} we prove that the answer on aforementioned question is
positive when the adaptive estimation of partial derivatives of a function belonging to an
anisotropic Nikol'skii class in $\bL_p(\bR^d)$ is considered  under an arbitrary $\bL_p$-loss and in an arbitrary dimension $d\geq 1$.

\paragraph{Acknowledgement.} The author is grateful to A. Goldenshluger and Yu. Golubev for fruitful discussions.

\section{Selection scheme for solving of {\it  Problem A}.}
\label{sec:ell-selector}

Let $\mH_n, n\in\bN^*$, be a sequence of countable subsets of $\mH$.
Let $\{\widehat{A}_h, \mh\in \mH\}$ and $\{\widehat{A}_{\mh,\eta},\; \mh,\eta\in \mH\}$ be the families of $X_1^{(n)}$-measurable $\mS_1$-valued mappings
possessing the  properties formulated below. Both $\widehat{A}_h$ and $\widehat{A}_{\mh,\eta}$  depend usually on $n$ but we will omit this dependence for the sake of simplicity of   notations. Let $\e_n\to 0, n\to\infty$ be a given sequence.


Suppose there exist  collections of $\mS_1$-valued  functionals  $\{\Lambda_\mh(f),\; \mh\in \mH\}$, $\{\Lambda_{\mh,\eta}(f),\; \mh,\eta\in \mH\}$ and a collection of {\it positive} $X_1^{(n)}$-measurable random variables $\Psi_n=\{\Psi_n(\mh),\; \mh\in \mH\}$  for which the following conditions hold. (The functionals  $\Lambda_\mh$ and $\Lambda_{\mh,\eta}$  may depend  on $n$ (not necessarily) but we will omit this dependence in the notations.)

\medskip

${\bf A^{\text{permute}}.}$\quad {\it $\widehat{A}_{\mh,\eta}\equiv\widehat{A}_{\eta,\mh}$,  for any $\mh, \eta\in\mH$.}

\medskip

${\bf A^{\text{upper}}.}$ {\it For any $n\geq 1$}
\begin{eqnarray*}
&&(\mathbf{i}) \;\quad\sup_{f\in\bF}\bE^{(n)}_{1,f}\bigg(\sup_{\mh\in \mH_n}\Big[\ell\big(\widehat{A}_\mh,\Lambda_\mh(f)\big)-\Psi_n(\mh)\Big]_+^{q}\bigg)\leq \e_n^{q};
\\*[1mm]
&&(\mathbf{ii}) \quad\sup_{f\in\bF}\bE^{(n)}_{1,f}\bigg(\sup_{\mh,\eta\in \mH_n}\Big[\ell\big(\widehat{A}_{\mh,\eta},\Lambda_{\mh,\eta}(f)\big)-\big\{\Psi_n(\mh)\wedge\Psi_n(\eta)\big\}\Big]_+^{q}\bigg)\leq \e_n^{q}.
\end{eqnarray*}

\begin{remark}
\label{r1}

Often the collection $\{\Psi_n(\mh),\; \mh\in \mH\}$ satisfying the hypothesis ${\bf A^{\text{upper}}}$ is not random. This is typically the case when a statistical problem is studied in the framework of
white gaussian noise or regression model. See also Lemma \ref{lem:extract-from-GL11} below, case $p\leq 2$.

\end{remark}

\subsection{Discussion.} For many statistical models and problems $\Lambda_\mh(f)=\bE^{(n)}_{1,f}\big(\widehat{A}_\mh\big)$ and $\Lambda_{\mh,\eta}(f)=\bE^{(n)}_{1,f}\big(\widehat{A}_{\mh,\eta}\big)$.  In this case
$
\ell\big(\widehat{A}_\mh,\Lambda_\mh(f)\big)$ and $\ell\big(\widehat{A}_{\mh,\eta},\Lambda_{\mh,\eta}(f)\big)$
can be viewed as stochastic errors related to the estimators $\widehat{A}_\mh$ and $\widehat{A}_{\mh,\eta}$ respectively. Hence, following the terminology used
in \cite{lepski2016} we can say that $\{\Psi_n(\mh),\; \mh\in \mH\}$ and $\{\Psi_n(\mh)\wedge\Psi_n(\eta),\; \mh,\eta\in \mH\}$ are upper functions of level $\e_n$ for the collection of corresponding stochastic errors.

The development of  the probabilistic tools allowing to control the behavior of stochastic errors related to statistical procedures is the subject of vast literature, see for instance the books  \cite{wellner}, \cite{van-de-Geer},
\cite{massart07}. The inequalities similar to those appeared in the hypothesis ${\bf A^{\text{upper}}}$  can be found in \cite{Ostrovskii} and \cite{lepski-a,lepski-b,lepski-c}.
The upper functions for $\bL_p$-norm of "kernel-type" empirical and gaussian processes  were obtained in  \cite{GL3} and \cite{lepski2016}.

We conclude that from the one hand the verification of the hypothesis ${\bf A^{\text{upper}}}$ is in some sense necessary task when the oracle approach or adaptive estimation is considered. On the other hand the very developed probabilistic machinery  is in our disposal.

Let us also remark that if the analogues of hypotheses ${\bf A^{\text{upper}}}$ and ${\bf A^{\text{permute}}}$ can be checked in {\it Problem B}, then the selection rule from the estimator's collection $\{\widehat{B}_h, \mh\in \mH\}$ presented below will provide the solution of {\it Problem B}.
However in some cases the verification of these hypotheses is much more difficult in {\it Problem B} than in {\it Problem A} and it is
one of the reasons why we propose to proceed differently.

\smallskip

Let us now discuss some examples of estimator's collections for which the hypothesis ${\bf A^{\text{permute}}}$ is fulfilled.

\begin{example}
\label{ex:1}
Let $\mathbf{A}\big(\cH^d\big)$ and $\mathbf{B}\big(\cH^d\big)$ be the families defined in (\ref{eq:def-M-est})
and (\ref{eq:def-K-est}). For any $\vec{h},\vec{h}^\prime\in\cH^d$ set
$\vec{h}\vee\vec{h}^\prime=(h_1\vee h^\prime_1,\ldots,h_d\vee h^\prime_d)$ and introduce
$$
\widehat{A}_{\vec{h},\vec{h}^\prime}(\cdot)=\widehat{A}_{\vec{h}\vee\vec{h}^\prime}(\cdot),\quad
\widehat{B}_{\vec{h},\vec{h}^\prime}(\cdot)=\widehat{B}_{\vec{h}\vee\vec{h}^\prime}(\cdot).
$$
It is obvious that the hypothesis ${\bf A^{\text{permute}}}$ is verified for these estimators. The selection rules
based on this construction of families of auxiliary estimators  can be found
in \cite{lepski-kerk}, \cite{lepski-kerk-08}, \cite{rebelles15}.

\end{example}

\begin{example}
\label{ex:2}
Consider either the density model generated by the observation $X^{(n)}=(X_1,\ldots,X_n)$, where $X_i\in\bR^d, i=1,\ldots n,$ are i.i.d. random vectors or the observation model (\ref{eq:observation-scheme}).

Let $K_{\vec{h}}(\cdot)$ and $M\big(\cdot, \vec{h}\big), \vec{h}\in\cH^d$, be defined in (\ref{eq:def-kernel-K}) and (\ref{eq:def-kernel-M})
 respectively. For any $\vec{h},\vec{h}^\prime\in\cH^d$ set
 $$
 K_{\vec{h},\vec{h}^\prime}(\cdot)=\int_{\bR^d}K_{\vec{h}}(\cdot-t)K_{\vec{h}^\prime}(t)\rd t,\quad M\big(\cdot, \vec{h},\vec{h}^\prime\big)=\int_{\bR^d}M\big(\cdot-t, \vec{h}\big)M\big(t, \vec{h}^\prime\big)\rd t
 $$
 and define
$$
\widehat{A}_{\vec{h},\vec{h}^\prime}(\cdot)=\sum_{i=1}^nK_{\vec{h},\vec{h}^\prime}(X_i-\cdot)\; \text{or}\;\;\widehat{A}_{\vec{h},\vec{h}^\prime}(\cdot)=\sum_{i=1}^nK_{\vec{h},\vec{h}^\prime}(Z_i-\cdot),\quad \widehat{B}_{\vec{h},\vec{h}^\prime}(\cdot)=\sum_{i=1}^nM\big(Z_i-\cdot, \vec{h},\vec{h}^\prime\big).
$$
Since, $K_{\vec{h},\vec{h}^\prime}\equiv K_{\vec{h}^\prime,\vec{h}}$ and $M\big(\vec{h},\vec{h}^\prime\big)\equiv M\big(\vec{h}^\prime,\vec{h}\big)$ the hypothesis ${\bf A^{\text{permute}}}$ holds. The selection rules based on this construction  of auxiliary estimators  were proposed
in \cite{LepLev98}, \cite{GL08}, \cite{GL09}, \cite{GL11}, \cite{comte}, \cite{patricia2014}, \cite{GL13}.

\end{example}

\begin{example}
\label{ex:3}
Let $\cD$ be a set endowed with the Borel measure $\mu$ and let
$\left\{\psi_{\blk}, \blk\in\bN^d\right\}$
 be an orthogonal basis in $\bL_2\big(\cD,\mu\big)$ possessing the following
properties.
$$\psi_{\blo}\equiv \mathbf{c}\neq 0,\quad \int_{\cD}\psi_{\blk} (t)\mu(\rd t)=0, \;\forall\blk\neq \blo=(0,\ldots 0).
$$
Let  $\bT=\{\tau: \tau=\big(\tau_{\blk},\blk\in\bN^d\big)\}$ be a given subset of $l_2$ such  that $\tau_{\blo}\mathbf{c}^2\mu(\cD)=1$ for all
$\tau\in\bT$. Introduce
$$
K_\tau(t,x)=\sum_{\blk\in\bN^d}\tau_{\blk}\psi_{\blk}(t)\psi_{\blk}(x),\;
\tau\in\bT,
$$
and define for any $\tau,\tau^{\prime}\in\bT$
$$
K_{\tau,\tau^\prime}(t,x)=\int_{\cD}K_\tau(t,y)K_\tau(y,x)\mu(\rd y).
$$
 Consider the statistical experiment  generated by the observation $X^{(n)}=(X_1,\ldots,X_n)$, where $X_i, i=1,\ldots n,$ are i.i.d. $\cD$-valued random variables having unknown  density $f$ with respect to the measure $\mu$. Introduce  the following collection of  estimators
$$
\widehat{A}_{\tau,\tau^\prime}(\cdot)=\sum_{i=1}^nK_{\tau,\tau^\prime}(X_i,\cdot),\quad \tau,\tau^\prime\in\bT.
$$
As it was shown in \cite{GL12}, $K_{\tau,\tau^\prime}\equiv K_{\tau^\prime,\tau}$ for any $\tau,\tau^\prime\in\bT$ and, therefore, the  hypothesis ${\bf A^{\text{permute}}}$ is fulfilled.

\end{example}

We remark that all estimator's construction discussed above can be applied in other statistical models where kernel-type estimators are used, e.g. white gaussian noise model and the regression model.

\subsection{$(\Psi_n,\ell)$-selection rule and corresponding oracle inequality.} Our objective is to propose the selection rule from an arbitrary collection $\mathbf{A}(\mH)=\{\widehat{A}_h, \mh\in \mH\}$ satisfying hypotheses  ${\bf A^{\text{permute}}}$ and  ${\bf A^{\text{upper}}}$, and establish for it the oracle inequality (\ref{eq:ell-oracle-inequality}).

Define for any $\mh\in \mH_n$
$$
\widehat{R}_n(\mh)=\sup_{\eta\in \mH_n}\Big[\ell\big(\widehat{A}_{\mh,\eta},\widehat{A}_\eta\big)- 2\Psi_n(\eta)\Big]_+.
$$
Let $\hat{\mh}^{(n)}\in\mH_n$ be an arbitrary $X_1^{(n)}-\text{measurable}$ random element satisfying
$$
\widehat{R}_n\big(\hat{\mh}^{(n)}\big)+2\Psi_n\big(\hat{\mh}^{(n)}\big)\leq\inf_{\mh\in \mH_n}\big\{\widehat{R}_n(\mh)+2\Psi_n(\mh)\big\}+\e^\prime_n,
$$
where $\e^\prime_n$ is an arbitrary sequence and we will assume that $\e^\prime_n\leq \e_n$.

Introduce the following notation: for any $f\in\bF$, $\mh\in\mH_n$ and  $n\geq 1$
$$
\cB^{(n)}_A(f,\mh)=\ell\big(\Lambda_\mh(f),A(f)\big)+2\sup_{\eta\in \mH_n}\ell\big(\Lambda_{\mh,\eta}(f),\Lambda_\eta(f)\big),\quad \psi_n(f,\mh)=\Big[\bE^{(n)}_{1,f}\big\{\Psi^q_n(\mh)\big\}\Big]^{\frac{1}{q}}.
$$

\begin{theorem}
\label{th1}
 Let  ${\bf A^{\text{permute}}}$ and ${\bf A^{\text{upper}}}$ be fulfilled. Then, for any $f\in\bF$ and   $n\geq 1$
$$
\cR_{A}\big[\widehat{A}_{\hat{\mh}^{(n)}}, f\big]\leq \inf_{\mh\in\mH_n}\left\{\cB^{(n)}_A(f,\mh)+5\psi_n(f,\mh)\right\}+6\e_n.
$$

\end{theorem}

Note  that the proposed  $(\Psi_n,\ell)$ selection rule is built using only observation part $X_1^{(n)}$. This  means that the application
of the result presented in Theorem \ref{th1} does not require any splitting of data since one can formally suppose that $X_1^{(n)}$ is the original data set.

\subsection{Proof of Theorem \ref{th1}.} For the simplicity of notations throughout the proof  we will write $\hat{\mh}$ instead of $\hat{\mh}^{(n)}$.

 Fix $\mh\in\mH_n$. We have obviously in view of the definition of $\widehat{R}_n(\cdot)$
\begin{eqnarray}
\label{eq1}
\ell\big(\widehat{A}_{\hat{\mh}},\widehat{A}_{\mh,\hat{\mh}}\big)&\leq& 2\Psi_n(\hat{\mh})+\Big[\ell\big(\widehat{A}_{\hat{\mh}},\widehat{A}_{\mh,\hat{\mh}}\big)-2\Psi_n(\hat{\mh})\Big]_+
\nonumber\\
&\leq& 2\Psi_n(\hat{\mh})+\widehat{R}_n(\mh).
\end{eqnarray}
Here we have also used that $\hat{\mh}\in\mH_n$ by its definition.
Taking into account  ${\bf A^{\text{permute}}}$ we get
\begin{eqnarray}
\label{eq2}
\ell\big(\widehat{A}_{\mh},\widehat{A}_{\mh,\hat{\mh}}\big)&=&\ell\big(\widehat{A}_{\mh},\widehat{A}_{\hat{\mh},\mh}\big)\leq 2\Psi_n(\mh)+\Big[\ell\big(\widehat{A}_{\mh},\widehat{A}_{\hat{\mh},\mh}\big)-2\Psi_n(\mh)\Big]_+
\nonumber\\
&\leq& 2\Psi_n(\mh)+\widehat{R}_n(\hat{\mh}).
\end{eqnarray}
We get from (\ref{eq1}),  (\ref{eq2}) and the definition of $\hat{\mh}$
\begin{eqnarray}
\label{eq3}
\ell\big(\widehat{A}_{\hat{\mh}},\widehat{A}_{\mh,\hat{\mh}}\big)+\ell\big(\widehat{A}_{\mh},\widehat{A}_{\mh,\hat{\mh}}\big)&\leq& \widehat{R}_n(\hat{\mh})+2\Psi_n(\hat{\mh})
+\widehat{R}_n(\mh)+2\Psi_n(\mh)
\nonumber\\
&\leq& 2\big\{\widehat{R}_n(\mh)+2\Psi_n(\mh)\big\}+\e^\prime_n.
\end{eqnarray}
We have in view of the triangle inequality for any $\mh\in \mH_n$
\begin{eqnarray}
\label{eq4}
\widehat{R}_n(\mh)\leq \sup_{\eta\in \mH_n}\ell\big(\Lambda_{\mh,\eta}(f),\Lambda_\eta(f)\big)+\xi_1+\xi_2,
\end{eqnarray}
where we have put
\begin{eqnarray}
\label{eq:def-xi_1-xi_2}
\xi_1=\sup_{\eta\in \mH_n}\Big[\ell\big(\widehat{A}_\eta,\Lambda_\eta\big)-\Psi_n(\eta)\Big]_+,\quad \xi_2=\sup_{\mh,\eta\in \mH_n}\Big[\ell\big(\widehat{A}_{\mh,\eta},\Lambda_{\mh,\eta}\big)-\big\{\Psi_n(\mh)\wedge\Psi_n(\eta)\big\}\Big]_+.&&
\end{eqnarray}
Thus, we obtain from  (\ref{eq3}) and (\ref{eq4}) for any $\mh\in \mH_n$
\begin{eqnarray}
\label{eq5}
\ell\big(\widehat{A}_{\hat{\mh}},\widehat{A}_{\mh,\hat{\mh}}\big)+\ell\big(\widehat{A}_{\mh},\widehat{A}_{\mh,\hat{\mh}}\big)\leq 2\sup_{\eta\in \mH_n}\ell\big(\Lambda_{\mh,\eta}(f),\Lambda_\eta(f)\big)+4\Psi_n(\mh)+2\xi_1+2\xi_2+\e^\prime_n.&&
\end{eqnarray}
Obviously for any $\mh\in \mH_n$
$$
\ell\big(\widehat{A}_{\mh},A(f)\big)\leq \ell\big(\Lambda_\mh(f),A(f)\big)+\Psi_n(\mh)+\xi_1.
$$
It yields together with (\ref{eq5}) by the triangle inequality
\begin{eqnarray}
\label{eq6}
\ell\big(\widehat{A}_{\hat{h}},A(f)\big)\leq \cB^{(n)}_A(f,\mh)+5\Psi_n(\mh)+3\xi_1+2\xi_2+\e^\prime_n,\quad\forall\mh\in\mH_n,\;\forall f\in\bF.
\end{eqnarray}
Taking into account the hypothesis ${\bf A^{\text{upper}}}$ we get for any $\mh\in\mH_n$ and any $f\in\bF$
$$
\bigg\{\bE^{(n)}_{1,f}\Big[\ell\big(\widehat{A}_{\hat{\mh}},A(f)\big)\Big]^q\bigg\}^{\frac{1}{q}}\leq \cB^{(n)}_A(f,\mh)+5\psi_n(f,\mh)+6\e_n.
$$
Here we have used that $\e^\prime_n\leq \e_n$.
Noting that the left hand side of the obtained inequality is independent of $\mh$ we come to the assertion of the theorem.
\epr

\section{"Plug-in" estimator for solving of {\it  Problem B} and its properties.}
\label{sec:rho-selector}

Let $\mathbf{B}(\mH)=\{\widehat{B}_\mh, \mh\in \mH\}$ be the family of $X_2^{(n)}$-measurable $\mS_2$-valued mappings.
The goal of this section is to bound from above the risk of "{\it plug-in}" estimator $\widehat{B}_{\hat{\mh}^{(n)}}$, where $\hat{\mh}^{(n)}$
is obtained by $(\Psi_n,\ell)$-selection rule.

We  add  the following assumption to the hypothesis ${\bf A^{\text{upper}}}$.

\smallskip

$\boldsymbol{\cA^{\text{upper}}.}$ {\it There exist the constant $C_{\Psi}\geq 1$ such that  for any $n\geq 1$}
\begin{equation*}
\sup_{f\in\bF}\bE^{(n)}_{1,f}\bigg(\sup_{\mh\in \mH_n}\Big[\Psi_n(\mh)-C_{\Psi}\psi_n(f,\mh)\Big]_+^{q}\bigg)\leq \e_n^{q},\quad
\sup_{f\in\bF}\bE^{(n)}_{1,f}\bigg(\sup_{\mh\in \mH_n}\Big[\psi_n(f,\mh)-C_{\Psi}\Psi_n(\mh)\Big]_+^{q}\bigg)\leq \e_n^{q}.
\end{equation*}
We note that the hypothesis $\boldsymbol{\cA^{\text{upper}}}$ is fulfilled with $C_{\Psi}=1$ is $\Psi_n(\mh)$ is non-random and in this case $\psi_n(\cdot,\cdot)$  is independent  of $f$.
Actually the hypothesis $\boldsymbol{\cA^{\text{upper}}}$ guarantees  that the random function $\Psi_n(\cdot)$ is well-concentrated around some non-random mapping. We would like to stress that for all known to the author problems in order  to check  ${\bf A^{\text{upper}}}$ one first verifies that the required inequalities hold for the "non-random" mapping $\psi_n(f,\mh)$. However it may depend on unknown $f$ and, therefore, cannot be used in the estimation construction. In this case this quantity is replaced by its
empirical counterpart $\Psi_n(\mh)$ satisfying the hypothesis $\boldsymbol{\cA^{\text{upper}}}$.

\smallskip

Let  $\{\Upsilon_\mh,\; \mh\in \mH\}$
be a collection of $\mS_2$-valued  functionals. Set for any $\mh\in\mH$,  $f\in\bF$ and $n\geq 1$
\begin{eqnarray}
\label{eq:def-calE}
\cE_n(f,\mh)=\Big(\bE^{(n)}_{2,f}\big\{\rho^{q}\big(\widehat{B}_{\mh},
\Upsilon_\mh(f)\big)\big\}\Big)^{\frac{1}{q}}.
\end{eqnarray}
Introduce  the following set of parameters.
\begin{eqnarray*}
\mV_n(f)&=&\Big\{\mh\in\mH_n:\; \psi_n(f,\mh)< 2C^2_\Psi\inf_{\eta\in\mH_n}\big[\cB^{(n)}_A(f,\eta)+2\psi_n(f,\eta)\big]\Big\};
\\*[0mm]
\mU_n(f)&=&\Big\{\mh\in\mH_n:\; \ell\big(\Lambda_{\mh},A(f)\big)< 4C_\Psi\inf_{\eta\in\mH_n}\big[\cB^{(n)}_A(f,\eta)+2\psi_n(f,\eta)\big]\Big\}.
\end{eqnarray*}
Set $\delta_n=\inf_{f\in\bF}\inf_{\mh\in\mH_n}\psi(f,\mh)$ and define for any $f\in\bF$ and $n\geq 1$
\begin{eqnarray*}
\tau_n(f)&=&\sup_{\mh\in\mV_n(f)}\cE_n(f,\mh)+\big(5\e_n/\delta_n)
\sup_{\mh\in\mH_n}\cE_n(f,\mh);
\\*[1mm]
\nu_n(f)&=&\sup_{\mh\in\mU_n(f)}\rho\big(\Upsilon_{\mh},B(f)\big)
+\big(30\e_n/\delta_n)\sup_{\mh\in\mH_n}\rho\big(\Upsilon_{\mh},B(f)\big).
\end{eqnarray*}

Let $\hat{\mh}^{(n)}$ is obtained by $(\Psi_n,\ell)$-selection rule with $\e_n^\prime\leq\delta_n/4$.

\begin{theorem}
\label{th2}
Let  ${\bf A^{\text{permute}}}$,  ${\bf A^{\text{upper}}}$ and  $\boldsymbol{\cA^{\text{upper}}}$  hold.
Then, for any $f\in\bF$ and   $n\geq 1$
$$
\cR_{B}\big[\widehat{B}_{\hat{\mh}^{(n)}}, f\big]\leq \nu_n(f)+\tau_n(f).
$$
\end{theorem}
We will see that in particular examples the sequence $\e_n$  decreases to zero very rapidly. It allows often to assert that
$$
\kappa_n:=\big(5\e_n/\delta_n)\sup_{f\in\bF}\bigg(
\sup_{\mh\in\mH_n}\cE_n(f,\mh)+6\sup_{\mh\in\mH_n}\rho\big(\Upsilon_{\mh},B(f)\big)
\bigg)\to 0,\;n\to\infty
$$
very fast and the statement of the theorem is reduced to
\begin{eqnarray}
\label{eq:cor-to-th2}
\cR_{B}\big[\widehat{B}_{\hat{\mh}^{(n)}}, f\big]\leq \sup_{\mh\in\mU_n(f)}\rho\big(\Upsilon_{\mh},B(f)\big)
+\sup_{\mh\in\mV_n(f)}\cE_n(f,\mh)+\kappa_n.
\end{eqnarray}
\begin{remark}
\label{rem4}
Note that the sets $\mV_n(f)$ and $\mU_n(f)$  are completely determined by the quantities  appeared in the solution of {\it Problem A} while $\cE_n(f,\mh)$ and
$\rho\big(\Upsilon_{\mh},B(f)\big)$ represent respectively the standard deviation of stochastic error  and approximation error of the estimator $\widehat{B}_\mh$ used in {\it Problem B}. The functionals  $\tau_n(f)$ and $\nu_n(f)$ explain how all these quantities are related.
\end{remark}

Although, (\ref{eq:cor-to-th2}) is simpler to analyze than the assertion of the theorem, even this oracle inequality is less natural than those proved in Theorem
\ref{th1}. Next result allows to understand better the properties of the estimator $\widehat{B}_{\hat{\mh}^{(n)}}$.
\begin{proposition}
\label{prop1}
Let  ${\bf A^{\text{permute}}}$ and ${\bf A^{\text{upper}}}$ hold. Suppose additionally that
there exists a constant $C_\ell>0$ such that for any $f\in\bF$,  $n\geq 1$ and $\mh\in\mH_n$
\begin{eqnarray}
\label{eq:hyp-cor1}
\rho\big(\Upsilon_\mh(f),B(f)\big)\leq C_\ell\;\ell\big(\Lambda_\mh(f),A(f)\big).
\end{eqnarray}
Assume also that there exists a constant $C_{\cE}>0$ such that for any   $n\geq 1$ and $\mh\in\mH_n$
\begin{eqnarray}
\label{eq:hyp-cor2}
\sup_{f\in\bF}\cE_n(f,\mh)\leq C_{\cE}\Psi_n(\mh).
\end{eqnarray}
Then, for any $f\in\bF$ and  $n\geq 1$
$$
\cR_{B}\big[\widehat{B}_{\hat{\mh}^{(n)}}, f\big]\leq \inf_{\mh\in\mH_n}\left\{(2C_\ell+C_{\cE})\cB^{(n)}_A(f,\mh)+(7C_\ell+2C_{\cE})\psi_n(f,\mh)\right\}+
(10C_\ell+3C_{\cE})\e_n.
$$

\end{proposition}

We remark that the assertions of Theorem \ref{th1}  and of Proposition \ref{prop1} differ only by the numerical constant. We would like to emphasize that although Proposition \ref{prop1} holds under quite restrictive assumption (\ref{eq:hyp-cor1}) it is useful for some problems studied in
Section \ref{sec:deconv}. Moreover, it does not require the hypothesis $\boldsymbol{\cA^{\text{upper}}}$.

\subsection{Proof of Theorem \ref{th2}.}

 For the simplicity of notations throughout the proof  we will write $\hat{\mh}$ instead of $\hat{\mh}^{(n)}$ and we break the proof on several steps.

$\mathbf{1^0}.$
Let us prove that for any $f\in\bF$  and  $n\geq 1$
\begin{eqnarray}
\label{eq:assert1}
 \bP^{(n)}_{1,f}\Big\{\hat{\mh}\notin\mU_n(f)\Big\}\leq (30\e_n/\delta_n)^q.
\end{eqnarray}
First note that in view of the triangle inequality
\begin{eqnarray}
\label{eq7}
 \ell\big(\Lambda_{\hat{\mh}},A(f)\big)&\leq& \ell\big(\widehat{A}_{\hat{\mh}},A(f)\big)+\ell\big(\widehat{A}_{\hat{\mh}},\Lambda_{\hat{\mh}}\big)
\nonumber \\
 &\leq& \ell\big(\widehat{A}_{\hat{\mh}},A(f)\big)+\Psi_n\big(\hat{\mh}\big)+
 \sup_{\eta\in\mH_n}\Big[\ell\big(\widehat{A}_{\eta},\Lambda_{\eta}\big)-\Psi_n(\eta)\Big]
 \nonumber \\
 &=& \ell\big(\widehat{A}_{\hat{\mh}},A(f)\big)+\Psi_n\big(\hat{\mh}\big)+\xi_1.
\end{eqnarray}
To get the last inequality we have used  that $\hat{\mh}\in\mH_n$.
Next, we get from the definition of $\hat{\mh}$ and (\ref{eq4})
\begin{eqnarray}
\label{eq8}
 \Psi_n\big(\hat{\mh}\big)
&\leq& \widehat{R}_n(\hat{\mh})+2\Psi_n(\hat{\mh})\leq
\widehat{R}_n(\mh)+2\Psi_n(\mh)+\e_n
\nonumber\\*[2mm]
&\leq&\sup_{\eta\in \mH_n}\ell\big(\Lambda_{\mh,\eta}(f),\Lambda_\eta(f)\big)+2\Psi_n(\mh)+\xi_1+\xi_2+\e^\prime_n
\nonumber\\*[1mm]
&\leq&\cB^{(n)}_A(f,\mh)+2\Psi_n(\mh)+\xi_1+\xi_2+\e^\prime_n.
\end{eqnarray}
Set $\xi_3=\sup_{\mh\in\mH_n}\big[\Psi_n(\mh)-C_{\Psi}\psi(f,\mh)\big]_+$. We deduce from (\ref{eq6}), (\ref{eq7}) and (\ref{eq8}) that
\begin{eqnarray}
\label{eq80}
 \ell\big(\Lambda_{\hat{\mh}},A(f)\big)
 &\leq& \inf_{\mh\in\mH}\big[2\cB^{(n)}_A(f,\mh)+7\Psi_n(\mh)\big]+5\xi_1+3\xi_2+2\e^\prime_n
 \nonumber\\
 &\leq& \inf_{\mh\in\mH}\big[2\cB^{(n)}_A(f,\mh)+7C_{\Psi}\psi(f,\mh)\big]+5\xi_1+3\xi_2+7\xi_3+2\e^\prime_n.
\end{eqnarray}
Hence for any $f\in\bF$ and   $n\geq 1$
\begin{eqnarray*}
\Big\{\hat{\mh}\notin\mU_n(f)\Big\}\subseteq \big\{
 5\xi_1+3\xi_2+7\xi_3+2\e^\prime_n\geq \delta_n\big\}\subseteq \big\{
 10\xi_1+6\xi_2+14\xi_3\geq \delta_n\big\}.
\end{eqnarray*}
To get the last inclusion we have used that $\e^\prime_n\leq\delta_n/4$. The statement  (\ref{eq:assert1}) follows now from the Markov inequality and the hypotheses ${\bf A^{\text{upper}}}$ and $\boldsymbol{\cA^{\text{upper}}}$.

\smallskip

$\mathbf{2^0}.$ For any $f\in\bF$ and  $n\geq 1$  the following is true:
\begin{eqnarray}
\label{eq:assert2}
 \bP^{(n)}_{1,f}\Big\{\hat{\mh}\notin\mV_n(f)\Big\}\leq (5\e_n/\delta_n)^q.
\end{eqnarray}
Indeed,  in view of (\ref{eq8})
$$
\Psi_n\big(\hat{\mh}\big)\leq 2\inf_{\mh\in\mH}\big[\cB^{(n)}_A(f,\mh)+C_{\Psi}\psi(f,\mh)\big]+\xi_1+\xi_2+2\xi_3+\e^\prime_n.
$$
Hence, putting $\xi_4=\sup_{\mh\in\mH_n}\big[\psi(f,\mh)-C_{\Psi}\Psi_n(\mh)\big]_+$ we obviously have
\begin{eqnarray*}
\psi(f,\hat{\mh})
\leq C_{\Psi}\Psi_n\big(\hat{\mh}\big)+\xi_4
\leq 2C_{\Psi}\inf_{\mh\in\mH}\big[\cB^{(n)}_A(f,\mh)+C_{\Psi}\psi(f,\mh)\big]
+C_{\Psi}\big[\xi_1+\xi_2+2\xi_3+\e^\prime_n\big]+\xi_4.
\end{eqnarray*}
Taking into account that $C_{\Psi}\geq 1$ and $\e^\prime_n\leq\delta_n/4$ we have
\begin{eqnarray*}
 \Big\{\hat{\mh}\notin\mV_n(f)\Big\}&=&\Big\{\psi(f,\hat{\mh})\geq 2C^2_\Psi\inf_{\eta\in\mH_n}\big[\cB^{(n)}_A(f,\eta)+2\psi_n(f,\eta)\big]\Big\}
 \nonumber\\
 &\subseteq& \Big\{\xi_1+\xi_2+2\xi_3++\xi_4+\e^\prime_n\geq 2\delta_n\Big\}\subseteq\Big\{\xi_1+\xi_2+2\xi_3+\xi_4\geq \delta_n\Big\}.
\end{eqnarray*}
The statement (\ref{eq:assert2}) follows now from the Markov inequality hypotheses ${\bf A^{\text{upper}}}$ and $\boldsymbol{\cA^{\text{upper}}}$.

\smallskip

$\mathbf{3^0}.$ Since $\hat{\mh}$ is $X_1^{(n)}$-measurable, $\{\widehat{B}_\mh, \mh\in\mH\}$ is $X_2^{(n)}$-measurable and $X_1^{(n)},X_2^{(n)}$ are independent
we get
\begin{eqnarray}
\label{eq11}
\Big(\bE^{(n)}_{f}\big\{\rho^q\big(\widehat{B}_{\hat{\mh}},\Upsilon_{\hat{\mh}}(f)\big)\big\}\Big)^{\frac{1}{q}}=
\Big(\bE^{(n)}_{1,f}\big\{\cE^q_n\big(f,\hat{\mh}\big)\big\}\Big)^{\frac{1}{q}}
\end{eqnarray}
Next, we obviously have
$$
\cE_n\big(f,\hat{\mh}\big)=\cE_n\big(f,\hat{\mh}\big)\mathrm{1}_{\hat{\mh}\in\mV_n(f)}+\cE_n\big(f,\hat{\mh}\big)\mathrm{1}_{\hat{\mh}\notin\mV_n(f)}\leq
\sup_{\mh\in\mV_n(f)}\cE_n\big(f,\mh\big)+\sup_{\mh\in\mH_n}\cE_n\big(f,\mh\big)\mathrm{1}_{\hat{\mh}\notin\mV_n(f)}.
$$
Hence, we get from (\ref{eq:assert2})
\begin{eqnarray}
\label{eq101}
\Big(\bE^{(n)}_{1,f}\big\{\cE^q_n\big(f,\hat{\mh}\big)\big\}\Big)^{\frac{1}{q}}\leq \sup_{\mh\in\mV_n(f)}\cE_n\big(f,\mh\big)+
\sup_{\mh\in\mH_n}\cE_n\big(f,\mh\big)\big(5\e_n/\delta_n)=\tau_n(f).
\end{eqnarray}

$\mathbf{4^0}.$ We have
\begin{eqnarray*}
\rho\big(\Upsilon_{\hat{\mh}}(f),B(f)\big)&=&\rho\big(\Upsilon_{\hat{\mh}}(f),B(f)\big)\mathrm{1}_{\hat{\mh}\in\mU(f)}
+\rho\big(\Upsilon_{\hat{\mh}}(f),B(f)\big)\mathrm{1}_{\hat{\mh}\notin \mU(f)}
\\*[2mm]
&\leq&\sup_{\mh\in\mU(f)}\rho\big(\Upsilon_{\mh}(f),B(f)\big)
+\sup_{\mh\in\mH_n}\rho\big(\Upsilon_{\mh}(f),B(f)\big)\mathrm{1}_{\hat{\mh}\notin \mU(f)}.
\end{eqnarray*}
Hence, in view of (\ref{eq:assert1})  we obtain
\begin{eqnarray}
\label{eq1011}
&&\Big(\bE^{(n)}_{f}\big\{\rho^q\big(\Upsilon_{\hat{\mh}}(f),B(f)\big)\big\}\Big)^{\frac{1}{q}}\leq \nu_n(f).
\end{eqnarray}
It remains to note that in view of triangle inequality
$$
\cR_{B}\big[\widehat{B}_{\hat{\mh}^{(n)}}, f\big]\leq \Big(\bE^{(n)}_{f}\big\{\rho^q\big(\Upsilon_{\hat{\mh}}(f),B(f)\big)\big\}\Big)^{\frac{1}{q}}+
\Big(\bE^{(n)}_{f}\big\{\rho^q\big(\widehat{B}_{\hat{\mh}},\Upsilon_{\hat{\mh}}(f)\big)\big\}\Big)^{\frac{1}{q}}
$$
and the assertion of the theorem follows from (\ref{eq101}) and (\ref{eq1011}).
\epr

\subsection{Proof of Proposition \ref{prop1}.}

We have in view of the assumption (\ref{eq:hyp-cor1})
 for any $f\in\bF$
\begin{eqnarray}
\label{eq9}
\rho\big(\widehat{B}_{\hat{\mh}},B(f)\big)&\leq& \rho\big(\Upsilon_{\hat{\mh}}(f),B(f)\big)+\rho\big(\widehat{B}_{\hat{\mh}},\Upsilon_{\hat{\mh}}(f)\big)
\nonumber\\*[2mm]
&\leq& C_\ell\;\ell\big(\Lambda_{\hat{\mh}}(f),A(f)\big)+\rho\big(\widehat{B}_{\hat{\mh}},\Upsilon_{\hat{\mh}}(f)\big)
\end{eqnarray}
Since by construction $\hat{\mh}$ is $X_1^{(n)}$-measurable and $X_1^{(n)}$ and $X_2^{(n)}$ are independent
\begin{equation*}
\Big(\bE^{(n)}_{f}\big\{\ell^q\big(\Lambda_{\hat{\mh}}(f),A(f)\big)\big\}\Big)^{\frac{1}{q}}=
\Big(\bE^{(n)}_{1,f}\big\{\ell^q\big(\Lambda_{\hat{\mh}}(f),A(f)\big)\big\}\Big)^{\frac{1}{q}}.
\end{equation*}
We deduce from  (\ref{eq80})  and the hypothesis ${\bf A^{\text{upper}}}$ for any $f\in\bF$
\begin{eqnarray}
\label{eq10}
\Big(\bE^{(n)}_{f}\big\{\ell^q\big(\Lambda_{\hat{\mh}}(f),A(f)\big)\big\}\Big)^{\frac{1}{q}}
\leq \inf_{\mh\in\mH_n}\left\{2\cB^{(n)}_A(f,\mh)+7\psi_n(f,\mh)\right\}+10\e_n.&&
\end{eqnarray}
For any $\mh\in\mH_n$ and $f\in\bF$ we have in view of (\ref{eq8}) and the assumption (\ref{eq:hyp-cor2}),
$$
\cE_n\big(f,\hat{\mh}\big)\leq  C_{\cE}\Psi_n\big(\hat{\mh}\big) \leq C_{\cE}\big[\cB^{(n)}_A(f,\mh)+2\Psi_n(\mh)+\xi_1+\xi_2+\e_n\big],
$$
and, therefore, we deduce from the hypothesis ${\bf A^{\text{upper}}}$  and (\ref{eq11})
$$
\Big(\bE^{(n)}_{f}\big\{\ell^q\big(\widehat{B}_{\hat{\mh}},\Upsilon_{\hat{\mh}}(f)\big)\big\}\Big)^{\frac{1}{q}}\leq
C_{\cE}\big[\cB^{(n)}_A(f,\mh)+2\psi_n(f,\mh)\big]+3C_{\cE}\e_n.
$$
Since the left hand side of the obtained inequality is independent of $\mh$ we get
\begin{eqnarray}
\label{eq13}
\Big(\bE^{(n)}_{f}\big\{\ell^q\big(\widehat{B}_{\hat{\mh}},\Upsilon_{\hat{\bmh}}(f)\big)\big\}\Big)^{\frac{1}{q}}\leq
C_{\cE}\inf_{\mh\in\mH_n}\left\{\cB^{(n)}_A(f,\mh)+2\psi_n(f,\mh)\right\}+3C_{\cE}\e_n.
\end{eqnarray}
The assertion of the theorem follows now from (\ref{eq9}), (\ref{eq10}), (\ref{eq13})  and the triangle inequality.

\epr

\subsection{\bf From oracle approach to adaptive estimation.}

In  this section we discuss how to use the oracle approach in adaptive estimation over the given scale of functional classes.
We present  several results concerning adaptation which will be directly deduced from Theorems \ref{th1}, \ref{th2} and Proposition \ref{prop1}.

Let $\big\{\bF_\ma,\ma\in\mA\big\}$ be a given collection of subsets of $\bF$ and suppose that an abstract oracle inequality (\ref{eq:ell-oracle-inequality}) is established. Define
$$
R_n\big(\bF_\ma\big)=\sup_{f\in\bF_\ma}\inf_{\mh\in\mH}\cA_{\ell}^{(n)}\left(f,\mh\right)+cr_n,\quad \ma\in\mA.
$$
We immediately deduce from (\ref{eq:ell-oracle-inequality}) that for any $\ma\in\mA$
$$
\limsup_{n\to\infty}R^{-1}_n\big(\bF_\ma\big)\sup_{f\in\bF_\ma}\cR_{A}\big[\widehat{A}_{\hat{\mh}}, f\big]\leq 1.
$$
Using modern statistical language we can state that the estimator $\widehat{A}_{\hat{\mh}}$ is $R_n$-adaptive, where $R_n=\{R_n\big(\bF_\ma\big), \ma\in\mA\}$ is the family of normalizations. If additionally one can prove that for any $\ma\in\mA$
$$
\liminf_{n\to\infty}R^{-1}_n\big(\bF_\ma\big)\inf_{\widetilde{Q}}\sup_{f\in\bF_\ma}\cR_{A}\big[\widetilde{Q}, f\big]>0,
$$
where infimum is taken over all $X^{(n)}$-measurable $\mS_1$-valued random mappings we can assert that the estimator $\widehat{A}_{\hat{\mh}}$ is {\it optimally adaptive} over the scale $\big\{\bF_\ma,\ma\in\mA\big\}$. The latter means that this estimator is simultaneously asymptotically minimax on each $\bF_\ma$.

\paragraph{Consequences of Theorem \ref{th1} and Proposition \ref{prop1}.}

Set $S_n=\{s_n(\bF_\ma), \ma\in\mA\}$, where
$$
s_n\big(\bF_\ma\big)=\sup_{f\in\bF_\ma}\inf_{\mh\in\mH_n}\big\{\cB^{(n)}_A(f,\mh)+\psi_n(f,\mh)\big\}
$$
and let  $\hat{\mh}^{(n)}$ is obtained by $(\Psi_n,\ell)$-selection rule.

The results below follow immediately  from Theorem \ref{th1} and Proposition \ref{prop1}.

\begin{theorem}
\label{th3}

Let  ${\bf A^{\text{permute}}}$ and  ${\bf A^{\text{upper}}}$  hold and assume  that
$\e_n=o\big(\inf_{\ma\in\mA}s_n(\bF_\ma)\big), n\to\infty$.

\smallskip

\noindent (1) Then, for any  $\ma\in\mA$
$$
\limsup_{n\to\infty}s^{-1}_n(\bF_\ma)\sup_{f\in\bF_\ma}\cR_{A}\big[\widehat{A}_{\hat{\mh}^{(n)}}, f\big]\leq 5.
$$
(2) If additionally  (\ref{eq:hyp-cor1}) and (\ref{eq:hyp-cor2}) are fulfilled, then for any  $\ma\in\mA$
$$
\limsup_{n\to\infty}s^{-1}_n(\bF_\ma)\sup_{f\in\bF_\ma}\cR_{B}\big[\widehat{B}_{\hat{\mh}^{(n)}}, f\big]\leq 7C_\ell+2C_{\cE}.
$$
\end{theorem}

Thus, the estimators $\widehat{A}_{\hat{\mh}^{(n)}}$ and $\widehat{B}_{\hat{\mh}^{(n)}}$ are simultaneously $S_n$-adaptive. Hence, if $\widehat{A}_{\hat{\mh}^{(n)}}$ is optimally adaptive in {\it Problem A}  one can expect that  $\widehat{B}_{\hat{\mh}^{(n)}}$ is optimally adaptive in {\it Problem B}.

\paragraph{Adaptive analogue of Theorem \ref{th2}.}

Define for any $\ma\in\mA$
\begin{equation}
\label{eq:def-psi_h(F_ma)}
\psi_n\big(\bF_\ma,\mh\big)=\inf_{f\in\bF_\ma}\psi_n(f,\mh)
\end{equation}
and introduce the following sets of parameters
\begin{equation*}
\cV_n(\ma)=\big\{\mh\in\mH_n:\; \psi_n(\bF_\ma,\mh)< 4C^2_{\Psi}s_n(\bF_\ma)\big\},\quad
\displaystyle{\cU_n(\ma,f)=\big\{\mh\in\mH_n:\; \ell\big(\Lambda_{\mh},A(f)\big)< 8C_{\Psi}s_n(\bF_\ma)\big\}}.
\end{equation*}
It is obvious that $\mV_n(f)\subseteq\cV_n(\ma)$ and $\mU_n(f)\subseteq\cU_n(\ma,f)$ for any $f\in\bF_\ma$. Set finally
\begin{eqnarray*}
\varphi_n(\bF_\ma)&=&\sup_{f\in\bF_\ma}\sup_{\mh\in\cV_n(\ma)}\cE_n(f,\mh)+
\sup_{f\in\bF_\ma}\sup_{\mh\in\cU_n(\ma,f)}\rho\big(\Upsilon_{\mh},B(f)\big),
\\
\kappa_n(\bF_\ma)&=&\big(5\e_n/\delta_n)\sup_{f\in\bF_\ma}\Big[
\sup_{\mh\in\mH_n}\cE_n(f,\mh)+6\sup_{\mh\in\mH_n}\rho\big(\Upsilon_{\mh},B(f)\big)
\Big]
\end{eqnarray*}
and let  $\hat{\mh}^{(n)}$ is obtained by $(\Psi_n,\ell)$-selection rule.

The following statement is the direct consequence of Theorem \ref{th2}.

\begin{theorem}
\label{th4}
Let  ${\bf A^{\text{permute}}}$,  ${\bf A^{\text{upper}}}$ and ${\boldsymbol \cA^{\text{upper}}}$ hold.

Assume that
$
\kappa_n(\bF_\ma)=o\big(\varphi_n(\bF_\ma)\big), n\to\infty,
$
for any $\ma\in\mA$.
Then, for any $\ma\in\mA$
$$
\limsup_{n\to\infty}\sup_{f\in\bF_\ma}\varphi^{-1}_n(\bF_\ma)\cR_{B}\big[\widehat{B}_{\hat{\bmh}^{(n)}}, f\big]\leq 1.
$$
\end{theorem}

\subsection{Some computations in the density model.}
\label{sec:subsec-problem-A-density}

Our objective now is to give  an example of statistical model and problem in which we are able to check
 the hypotheses ${\bf A^{\text{permute}}}$,  ${\bf A^{\text{upper}}}$ and ${\boldsymbol \cA^{\text{upper}}}$.
We also gives some bounds for  the quantities $\cB^{(n)}_A(f,\mh)$ and $\psi_n(f,\mh)$ involved in the $(\Psi_n,\ell)$-oracle inequality.
It allows us, in particular, to compute the rate
$$
s_n\big(\bF_\ma\big)=\sup_{f\in\bF_\ma}\inf_{\mh\in\mH_n}\big\{\cB^{(n)}_A(f,\mh)+\psi_n(f,\mh)\big\}
$$
appeared  in all adaptive results in the case when $\bF_\ma=\bN_{p,d}\big(\vec{\beta},L\big)$, $\ma=\big(\vec{\beta},L\big)$, where $\bN_{p,d}\big(\vec{\beta},L\big)$ anisotropic Nikol'skii class.
The results presented in this section form basic tools  for the solution of problem studied in  Sections \ref{sec:deconv} and
\ref{sec:derivatives}.

 \smallskip

Let $T^{(m)}=(T_1,\ldots,T_m), m\in\bN^*,$ be i.i.d. random vectors taking values in $\bR^d$ and having density $\mathfrak{p}$ with respect to Lebesgue measure.
As before  $\bP^{(m)}_\mathfrak{p}$ and $\bE^{(m)}_\mathfrak{p}$ denote the probability law and mathematical expectation of  $T^{(m)}$.
The goal is to estimate $A(\mmp)=\mmp$  and  the quality of an estimator procedures is measured by $\bL_p$-risk, $1\leq p<\infty$, that is
$$
\cR_{A}\big[\widetilde{A}_m, \mmp\big]=\Big(\bE^{(m)}_\mmp\Big[\big\|\widetilde{A}_m-\mmp\big\|^p_p\Big]\Big)^{1/p},\quad \mmp\in\mP
$$
where $\mP$ is the set of all probability densities and
$
\ell^p(g)=\|g\|^p_p=\int_{\bR^d}|g(x)|^p\rd x.
$

\smallskip

Remind that $\cH^d=\big\{\big((2^{-k_1},\ldots,2^{-k_1}\big)), (k_1,\ldots,k_d)\in\bN^d\big\}$ is the diadic grid in $(0,1]^d$
and $K_{\vec{h}}(\cdot), \vec{h}\in\cH^d$ be defined in (\ref{eq:def-kernel-K}).
Introduce the collection of standard kernel estimators
\begin{equation}
\label{eq2:def-K-est-new}
\mathbf{A}\big(\cH_m^d\big)=\bigg\{\widehat{A}_{\vec{h}}(\cdot)=m^{-1}\sum_{i=1}^{m}K_{\vec{h}}(T_i-\cdot),\;\vec{h}\in\cH_m^d \bigg\},
\end{equation}
where
$
\cH_m^d=\big\{\vec{h}\in\cH^d: \ln(m)/m\leq V_{\vec{h}}\leq e^{-\sqrt{\ln(m)}}\big\}
$
and $V_{\vec{h}}=\prod_{j=1}^dh_j$.

\smallskip

From now on we will assume that the kernel $K$ involved in the definition of $K_{\vec{h}}$ is continuous, bounded, symmetric  function belonging to $\bL_1(\bR^d)$
and $\int K=1$.

\subsubsection{Verification of the hypotheses ${\bf A^{\text{permute}}}$,  ${\bf A^{\text{upper}}}$ and ${\boldsymbol \cA^{\text{upper}}}$.}
 \label{sec:subsec-verification-A}
 Recall that "$\star$" denotes the convolution operator and  introduce for any $\vec{h},\vec{\eta}\in\cH^d$
\begin{equation}
\label{eq:def-K-star-K}
\widehat{A}_{\vec{h},\vec{\eta}}(\cdot)=m^{-1}\sum_{i=1}^{m}\big[K_{\vec{h}}\star K_{\vec{\eta}}\big](T_i-\cdot).
\end{equation}
Set also for any $\vec{h}\in\cH^d$ and  $m\in\bN^*$
$$
\Delta_m\big(\vec{h}\big)=
\left\{
\begin{array}{lll}
&128\|K\|_1\|K\|_p\big(mV_{\vec{h}}\big)^{1/p-1},&\quad 1\leq p<2;
\\*[2mm]
&9\|K\|_1\|K\|_2\big(mV_{\vec{h}}\big)^{-1/2},&\quad p=2;
\end{array}
\right.
$$
and  if $p>2$
$$
\Delta_m\big(\vec{h}\big)=\frac{480p\|K\|_1}{\ln(p)}\bigg\{\frac{1}{\sqrt{m}}\bigg(\int_{\bR^d}\bigg[\frac{1}{m}\sum_{i=1}^{m}K^2_{\vec{h}}(Z_i-t)\bigg]^{p/2}\rd t\bigg)^{1/p}+2\|K\|_p\big(mV_{\vec{h}}\big)^{-1/2}\bigg\}.
$$
We remark that $\Delta_m(\cdot)$ is not random if $p\leq 2$.
It is obvious that the hypothesis ${\bf A^{\text{permute}}}$ is fulfilled.  Define the empirical process
$$
\xi_m\big(x,\vec{h}\big)=\frac{1}{m}\sum_{i=1}^{m}\Big[K_{\vec{h}}(T_i-x)-\bE^{(m)}_\mathfrak{p}K_{\vec{h}}(T_i-x)\Big], \quad x\in\bR^d.
$$

\noindent {\bf Auxiliary results.} The following statements are the consequences of the results established in Lemmas 1 and 2 in
\cite{GL11} and Theorem 3 in \cite{GL3}.
Set $\mP(D)=\{\mmp\in\mP: \|\mmp\|_\infty\leq D\}$, $D>0$.

\begin{lemma}
\label{lem:extract-from-GL11}
For any $1\leq p<\infty$ and for all $m\in\bN^*$
\begin{gather}
\label{eq1:lem1}
\sup_{\mmp\in\mP}\bE^{(m)}_\mathfrak{p}\bigg(\sup_{\vec{h}\in \cH^d_{m}}\Big[\big\|\xi_m\big(\cdot,\vec{h}\big)\big\|_p-\|K\|_1^{-1}
\Delta_m(\vec{h})\Big]_+^{q}\bigg)\leq \e^q_m(p),\quad p\in [1,2);
\\*[2mm]
\sup_{\mmp\in\mP(D)}\bE^{(m)}_\mathfrak{p}\bigg(\sup_{\vec{h}\in \cH^d_{m}}\Big[\big\|\xi_m\big(\cdot,\vec{h}\big)\big\|_p-\|K\|_1^{-1}
\Delta_m(\vec{h})\Big]_+^{q}\bigg)\leq \e^q_m(p),\quad p\in [2,\infty),
\label{eq2:lem1}
\end{gather}
where for any $p\geq 1$ and $D>0$
\begin{eqnarray}
\label{eq3:lem1}
\limsup_{m\to\infty}m^{a}\e_m(p)=0,\quad\forall a>0.
\end{eqnarray}
Moreover for any $p>2$, there exists $\mathfrak{D}>0$ completely determined by $K,p,D$ and $d$ such that
\begin{eqnarray}
\label{eq4-new:lem1}
\sup_{\mmp\in\mP(D)}\bE^{(m)}_\mathfrak{p}\bigg(\sup_{\vec{h}\in \cH^d_{m}}\Big[\Delta_m(\vec{h})-\mD\big(mV_{\vec{h}}\big)^{-1/2}\Big]_+^{q}\bigg)\leq \e^q_m(p).
\end{eqnarray}

\end{lemma}

Let us remark the following simple consequence of (\ref{eq3:lem1}) and (\ref{eq4-new:lem1}).
\begin{eqnarray}
\label{eq4:lem1}
\sup_{\mmp\in\mP(D)}\Big(\bE^{(m)}_\mathfrak{p}\Delta^{q}_m(\vec{h})\Big)^{1/q}\leq \mathfrak{D}
\big(mV_{\vec{h}}\big)^{-1/2},\quad\forall\vec{h}\in \cH^d_{m},\;\forall m\in\bN^*.
\end{eqnarray}
In order to get this inequality it suffices to note that $\big(mV_{\vec{h}}\big)^{-1/2}\geq m^{-1/2}$ for all $\vec{h}\in (0,1]^d$.

\paragraph{Verification of the hypothesis ${\bf A^{\text{upper}}}$.}

Set $\Lambda_{\vec{h}}(\mmp,\cdot)=\bE^{(m)}_{\mmp}\big\{\widehat{A}_{\vec{h}}(\cdot)\big\}$, $\vec{h}\in \cH^d_{m}$ and note that
\begin{equation}
\label{eq:Lambda-new}
\Lambda_{\vec{h}}(\mmp,\cdot)=\bE^{(m)}_\mathfrak{p}K_{\vec{h}}(T_1-\cdot)=\int_{\bR^d}K_{\vec{h}}(z-\cdot)\mmp(z)\rd z.
\end{equation}
In view of the assertions (\ref{eq1:lem1}) and(\ref{eq2:lem1}) of Lemma \ref{lem:extract-from-GL11}  we can conclude that
the hypothesis ${\bf A^{\text{upper}} (i)}$ is fulfilled for any $p\geq 1$.
In order to verify ${\bf A^{\text{upper}} (ii)}$ note first the following obvious identities.
\begin{eqnarray}
\label{eq:identities}
\widehat{A}_{\vec{h},\vec{\eta}}\equiv K_{\vec{h}}\star\widehat{A}_{\vec{\eta}}\equiv K_{\vec{\eta}}\star\widehat{A}_{\vec{h}},\quad\forall\vec{h},\vec{\eta}\in\cH^d.
\end{eqnarray}
Hence, putting $\Lambda_{\vec{h},\vec{\eta}}(\mmp,\cdot)=\bE^{(m)}_{\mmp}\Big\{\widehat{A}_{\vec{h},\vec{\eta}}(\cdot)\Big\},\;$
$\zeta^{(m)}_{\vec{h},\vec{\eta}}(\cdot)=\widehat{A}_{\vec{h},\vec{\eta}}(\cdot)-\Lambda_{\vec{h},\vec{\eta}}(\mmp,\cdot)$
and $\zeta^{(m)}_{\vec{h}}(\cdot)=\widehat{A}_{\vec{h}}(\cdot)-\Lambda_{\vec{h}}(\mmp,\cdot)$ we deduce from the Young inequality, \cite{folland}, Theorem 6.18,
$$
\Big\|\zeta^{(m)}_{\vec{h},\vec{\eta}}\Big\|_p\leq \|K\|_1\Big\|\zeta^{(m)}_{\vec{\eta}}\Big\|_p,\quad \Big\|\zeta^{(m)}_{\vec{h},\vec{\eta}}\Big\|_p
\leq \|K\|_1\Big\|\zeta^{(m)}_{\vec{h}}\Big\|_p\;\Rightarrow\;
\Big\|\zeta^{(m)}_{\vec{h},\vec{\eta}}\Big\|_p\leq \|K\|_1\Big(\Big\|\zeta^{(m)}_{\vec{\eta}}\Big\|_p\wedge\Big\|\zeta^{(m)}_{\vec{h}}\Big\|_p\Big).
$$
It yields for all $\vec{h},\vec{\eta}\in\cH_m^d$ in view of the obvious inequality $a\wedge c\leq b\wedge d+(a-b)_+\vee (c-d)_+$
$$
\Big(\Big\|\zeta^{(m)}_{\vec{h},\vec{\eta}}\Big\|_p-\Delta_{m}(\vec{h})\wedge\Delta_{m}(\vec{\eta})\Big)_+ \leq \|K\|_1
\sup_{\vec{h}\in \cH^d_{m}}\Big(\Big\|\zeta^{(m)}_{\vec{h}}\Big\|_p-\|K\|^{-1}_1\Delta_{m}(\vec{h})\Big)_+
$$
In view of the assertions (\ref{eq1:lem1}) and (\ref{eq2:lem1}) of Lemma (\ref{lem:extract-from-GL11}) we can conclude that
the hypothesis ${\bf A^{\text{upper}} (ii)}$ is fulfilled for any $p\geq 1$.

\paragraph{Verification of the hypothesis ${\boldsymbol \cA^{\text{upper}}}$.}

First, we note that $\Delta_m$ is not random if $p\in[1,2]$ and, therefore, in this case ${\boldsymbol \cA^{\text{upper}}}$ is checked with $C_\Psi=1$.

\smallskip

Next, if $p>2$ we have in view of the definition of $\Delta_m$ and (\ref{eq4:lem1}) for any $\mmp\in\mP(D)$
\begin{eqnarray}
\label{eq100}
\frac{960p\|K\|_1\|K\|_p}{\ln(p)}\big(mV_{\vec{h}}\big)^{-1/2}\leq\Big(\bE^{(m)}_\mathfrak{p}\Delta^{q}_m(\vec{h})\Big)^{1/q}\leq \mathfrak{D}
\big(mV_{\vec{h}}\big)^{-1/2},\quad\forall\vec{h}\in \cH^d_{m}.
\end{eqnarray}
Moreover, $\Delta_m\geq \frac{960p\|K\|_1\|K\|_p}{\ln(p)}\big(mV_{\vec{h}}\big)^{-1/2}$ and, therefore, ${\boldsymbol \cA^{\text{upper}}}$ is fulfilled in view of
the assertion (\ref{eq2:lem1}) of Lemma (\ref{lem:extract-from-GL11}) with
$
C_\Psi=\mD\vee\frac{960p\|K\|_1\|K\|_p}{\ln(p)}.
$

\subsubsection{Some bounds for the quantities $\cB^{(n)}_A(f,\mh)$ and $\psi_n(f,\mh)$.}
\label{sec:subsec-consec-th1}

Let us find upper estimate for  $\cB^{(n)}_A(f,\mh)$ and lower and upper estimates for $\psi_n(f,\mh)$, which both appear in $(\Psi_n,\ell)$-oracle inequality.
In the considered case  $f=\mmp,  n=m$,
$\mh=\vec{h}$, $A(f)=A(\mmp)=\mmp$, $\ell(\cdot)=\|\cdot\|_p$ and  $\Psi_n=\Delta_{m}$.

\smallskip

Note first that  (\ref{eq:identities}) implies that for any $x\in\bR^d$
\begin{eqnarray*}
\Lambda_{\vec{h},\vec{\eta}}(\mmp,x)-\Lambda_{\vec{\eta}}(\mmp,x)&=&K_{\vec{\eta}}\star\Lambda_{\vec{h}}(\mmp,x)-\Lambda_{\vec{\eta}}(\mmp,x)
\\
&=&
\int_{\bR_d}K_{\vec{\eta}}(x-y)\bigg(\int_{\bR^d}K_{\vec{h}}(z-y)\mmp(z)\rd z-\mmp(y)\bigg)\rd y.
\end{eqnarray*}
Applying the Young inequality we obtain for all $p\geq 1$
\begin{eqnarray}
\label{eq:20}
\cB^{(m)}_A\big(\mmp,\vec{h}\big)\leq (1+2\|K\|_1)\big\|K_{\vec{h}}\star \mmp- \mmp\big\|_p,\quad \forall \vec{h}\in\cH^d.
\end{eqnarray}
In view of assertion (\ref{eq4:lem1}) of Lemma \ref{lem:extract-from-GL11} and the definition of $\Delta_{m}$ we get
\begin{eqnarray}
\label{eq:21}
\underline{\mathbf{C}}\big(nV_{\vec{h}}\big)^{\frac{1}{p\wedge2}-1}\leq\psi_n\big(\mmp,\vec{h}\big):=
\Big[\bE^{(m)}_{\mmp}\big\{\Delta^p_{m}(\mh)\big\}\Big]^{\frac{1}{p}}\leq
\overline{\mathbf{C}}\big(nV_{\vec{h}}\big)^{\frac{1}{p\wedge2}-1},\quad\forall \mmp\in\mP_p,
\end{eqnarray}
where $\mP_p=\mP$ if $p<2$,  $\mP_p=\{\mmp\in\mP: \mmp\in\mP(D)\}$ if $p\geq 2$ and
$\underline{\mathbf{C}}=\overline{\mathbf{C}}=128\|K\|_1\|K\|_p$ if $p\in [1,2)$,
$\underline{\mathbf{C}}=\overline{\mathbf{C}}=9\|K\|_1\|K\|_2$ if $p=2$ and $\underline{\mathbf{C}}=\frac{960p\|K\|_1\|K\|_p}{\ln(p)}$,
$\overline{\mathbf{C}}=\mD$ if $p>2$.
We remark that the obtained lower and upper estimates are independent of $\mmp$.

\subsubsection{Computation of the rate on  anisotropic Nikol'skii class.}
\label{sec:subsec-consec-th3}

In view of   (\ref{eq:20}) and (\ref{eq:21}) we should bound from above
$$
\sup_{\mmp\in\bF_\ma}\inf_{\vec{h}\in H^d_m}\left\{\big\|K_{\vec{h}}\star \mmp- \mmp\big\|_p+\big(mV_{\vec{h}}\big)^{\frac{1}{p\wedge2}-1}\right\}.
$$
This quantity with $\bF_\ma=\bN_{p,d}\big(\vec{\beta},L\big),\;\ma=(\vec{\beta},L\big)$, where $\bN_{p,d}\big(\vec{\beta},L\big)$ is anisotropic Nikol'skii class was
evaluated in \cite{GL11}.

Let $(\mathbf{e}_1,\ldots,\mathbf{e}_d)$ denote the canonical basis of $\bR^d$. For function $G:\bR^d\to \bR^1$ and
real number $u\in \bR$ define
 the first order difference operator with step size $u$ in direction of the variable
$x_j$ by
$$
 \Delta_{u,j}G (x)=G(x+u\mathbf{e}_j)-G(x),\;\;\;j=1,\ldots,d.
$$
By induction,
the $k$-th order difference operator with step size $u$ in direction of the variable $x_j$ is
defined~as
\begin{equation*}
\label{eq:Delta}
 \Delta_{u,j}^kG(x)= \Delta_{u,j} \Delta_{u,j}^{k-1} G(x) = \sum_{l=1}^k (-1)^{l+k}\binom{k}{l}\Delta_{ul,j}G(x).
\end{equation*}
\begin{definition}
\label{def:nikolskii}
For given   $\vec{\beta}=(\beta_1,\ldots,\beta_d)\in (0,\infty)^d, L>0$ and $p\geq 1$
  we
say that function $G:\bR^d\to \bR^1$ belongs to the anisotropic
Nikolskii class $\bN_{p,d}\big(\vec{\beta},L\big)$ if  $\|G\|_{p}\leq L$ and for every $j=1,\ldots,d$ there exists natural number  $k_j>\beta_j$ such that
\begin{equation*}
\label{eq:Nikolski}
 \Big\|\Delta_{u,j}^{k_j} G\Big\|_{p} \leq L |u|^{\beta_j},\;\;\;\;
\forall u\in \bR,\;\;\;\forall j=1,\ldots, d.
\end{equation*}

\end{definition}

We will use the following specific kernel $K$ in the definition of the family $\mathbf{A}\big(\cH_m^d\big)$
[see, e.g., \cite{lepski-kerk} or \cite{GL11}].
\par
 Let  $\ms\in\bN^*$ be  fixed and let $w:\bR^1\to\bR^1$ satisfy   $\int w(y)\rd y=1$,
and $w\in\bC^{(\ms)}(\bR^1)$. Put
\begin{equation}
\label{eq:w-function}
 \cK_\ms(y)=\sum_{i=1}^\ms \binom{\ms}{i} (-1)^{i+1}\frac{1}{i}w\Big(\frac{y}{i}\Big),\qquad
 K(x)=\prod_{j=1}^d \cK_\ms(x_j),\;\;\;\;x=(x_1,\ldots,x_d).
\end{equation}


Let $\ms>0$ be an arbitrary but a priory chosen number and let the kernel $K$ is constructed in accordance with (\ref{eq:w-function}).
Next result can be found for instance in \cite{GL11}.

For any $\vec{\beta}\in (0,\ms]^d$, $L>0$ and $p\geq 1$ there exists $\boldsymbol{\upsilon}>0$ independent of $L$ such that
\begin{eqnarray}
\label{eq:23-adaptive-new}
\sup_{\mmp\in\bN_{p,d}(\vec{\beta},L)\cap\mP_p}\;\inf_{\vec{h}\in H^d_m}\left\{\big\|K_{\vec{h}}\star \mmp- \mmp\big\|_p+\big(mV_{\vec{h}}\big)^{\frac{1}{p\wedge2}-1}\right\}&\leq & \boldsymbol{\upsilon}\;s_m\big(\bN_{p,d}\big(\vec{\beta},L\big)\cap\mP_p\big)
\\
\nonumber
&=:&\boldsymbol{\upsilon}\;L^{\frac{1-1/(p\wedge 2)}{\beta+1-1/(p\wedge 2)}}
\;m^{-\frac{1-1/(p\wedge 2)}{1+(1/\beta)(1-1/(p\wedge 2)}}.
\end{eqnarray}
It is well-known that $s_m\big(\bN_{p,d}\big(\vec{\beta},L\big)\cap\mP_p\big)$ is the minimax rate of convergence on
$\bN_{p,d}\big(\vec{\beta},L\big)\cap\mP_p$.


\smallskip

\noindent{\bf Concluding remarks.} $(\Psi_n,\ell)$-selection rule with $\Psi_n=\Delta_m, n=m,$ and $\ell(\cdot)=\|\cdot\|_p$ is a particular case of the
selection scheme proposed in \cite{GL11}. Selected in accordance with this rule estimator from the collection (\ref{eq2:def-K-est-new}) is, in view of
(\ref{eq:20}), (\ref{eq:21}),  (\ref{eq:23-adaptive-new}) and the first assertion of Theorem \ref{th3}, optimally-adaptive over the collection of anisotropic Nikol'skii
classes. For the first time it was proved in \cite{GL11}, Theorem 4.

\section{Adaptive estimation in the generalized deconvolution model.}
\label{sec:deconv}

 Consider the following observation scheme introduced in \cite{LW15}:
\begin{equation}
\label{eq:observation-scheme-1}
Z_i=X_i+\epsilon_iY_i,\quad i=1,\ldots,n,
\end{equation}
where $X_i, i=1,\ldots,n,$ are {\it i.i.d.} $d$-dimensional  random vectors with common density $f$ to be estimated.
The noise variables  $Y_i, i=1,\ldots,n,$
are {\it i.i.d.} $d$-dimensional random vectors with known common density $g$. At last $\e_i\in\{0,1\}, i=1,\ldots,n,$
are {\it i.i.d.} Bernoulli random variables with $\bP(\e_1=1)=\alpha$, where $\alpha\in [0,1)$ is supposed to be known.
The sequences  $\{X_i, i=1,\ldots,n\}$, $\{Y_i, i=1,\ldots,n\}$ and $\{\epsilon_i, i=1,\ldots,n\}$ are supposed to be mutually independent.

Let us note that the case $\alpha=1$ corresponds to the pure deconvolution model  $Z_i=X_i+Y_i,\; i=1,\ldots,n$ discussed in Introduction, whereas the case $\alpha=0$ corresponds to the direct observation scheme $Z_i=X_i,\; i=1,\ldots,n$. The intermediate case $\alpha\in (0,1)$ can be treated as the mathematical modeling of the following situation. One part of the data, namely $(1-\alpha) n$, is observed without noise. If the indexes corresponding to these observations were known, the density $f$ could be estimated using only this part of the data, with the accuracy corresponding to the direct case. The main question we will address in this intermediate case is whether the same accuracy would be achievable if the latter information is not available?

\smallskip

Thus,
$
\rho^p(g)=\|g\|^p_p=\int_{\bR^d}|g(x)|^p\rd x,
$
$\bF=\mP$ is the set of all probability densities,  $B(f)=f$ and the quality of an estimator procedures is measured by $\bL_p$-risk
$$
\cR_{B}\big[\widetilde{B}_n, f\big]=\Big(\bE^{(n)}_f\Big[\big\|\widetilde{B}_n-f\big\|^p_p\Big]\Big)^{1/p},\quad f\in\mP.
$$
At last, let $X^{(n)}_1=(X_1,\ldots, X_{[n/2]})$ and $X^{(n)}_2=(X_{[n/2]+1},\ldots,X_n)$, where $[a]$ denotes the integer part of $a\in\bR$.

All  results presented in this section
will be a established under the following condition imposed on the distribution of the noise variable $Y_1$. In what follows for any $Q\in\bL_1(\bR^d)$ its Fourier  transform will be denoted by $\check{Q}$.

\begin{assumption}
\label{ass1:ass-on-noise-upper-bound}

1) if $\alpha\in (0,1)$ then there exists  $\varpi>0$ such that
$$
\big|1 - \alpha +\alpha\check{g}(t)\big|\geq \varpi, \quad \forall t\in\bR^d;
$$

2) if $\alpha=1$ then there exists  $\vec{\mu}=(\mu_1,\ldots,\mu_d)\in (0,\infty)^{d}$ and $G_1,G_2>0$  such that
\begin{eqnarray*}
 G_1\prod_{j=1}^d(1+t^2_j)^{-\frac{\mu_j}{2}}\leq |\check{g}(t)|\leq G_2\prod_{j=1}^d(1+t^2_j)^{-\frac{\mu_j}{2}},\quad\forall t\in\bR^d.&&
\end{eqnarray*}

\end{assumption}

We remark that Assumption \ref{ass1:ass-on-noise-upper-bound} is very weak  when $\alpha\in(0,1)$.
It is verified for many distributions, including centered  multivariate Laplace and Gaussian ones.
Note also that Assumption \ref{ass1:ass-on-noise-upper-bound} always holds  with $\varpi=1-2\alpha$ if $\alpha<1/2$.
Additionally, it holds  with $\varpi=1-\alpha$ if $\check{g}$ is a real positive function. The latter is true, in particular, for any probability law obtained by the even number of convolutions of a symmetric distribution with itself. As to the case $\alpha=1$ Assumption \ref{ass1:ass-on-noise-upper-bound} is s well-known in the literature condition
referred to {\it moderately  ill-posed}   statistical problem. In particular, it is checked for   multivariate Laplace and Gamma laws.

\smallskip

Below we introduce  families of estimators whose construction involved  kernel $K$.  Throughout this section without further mentioning we will additionally
assume that  $\check{K}$ is compactly supported and
$
\check{K}(t)=1
$
for all  $t\in[-1,1]^d$.

From now on $\mathbf{c_1},\mathbf{c_2},\ldots,$ denote the constants that may depends only on $K,p,D,\mathfrak{D},d$ and
 $\alpha,\varpi, G_1, G_2, \vec{\mu}$ appeared in Assumption \ref{ass1:ass-on-noise-upper-bound}. In particular they are independent of $f$ and $n$.


\subsection{Idea of estimator construction.}
\label{sec:subsec-idea}

Let  $K_{\vec{h}}(\cdot), \vec{h}\in (0,1]^d$, be defined in (\ref{eq:def-kernel-K})
and let  $M\big(\cdot,\vec{h}\big):\bR^d\to\bR$
satisfy the operator equation
\begin{eqnarray}
\label{eq2:def-kernel-M}
K_{\vec{h}}(y)=(1-\alpha)M\big(y,\vec{h}\big)+\alpha\int_{\bR^d}g(t-y)M\big(t,\vec{h}\big)\rd t,\quad  y\in\bR^d.&&
\end{eqnarray}
Recall that $V_{\vec{h}}=\prod_{j=1}^dh_j$,
$
\cH_m^d=\big\{\vec{h}\in\cH^d: \ln(m)/m\leq V_{\vec{h}}\leq e^{-\sqrt{\ln(m)}}\big\}, m\geq 1,
$
and introduce the following estimator's collection built from $X_2^{(n)}$
\begin{equation}
\label{eq2:def-M-est}
\mathbf{B}\big(\cH_{[n/2]}^d\big)=\bigg\{\widehat{B}_{\vec{h}}(\cdot)=(n-[n/2])^{-1}\sum_{i=[n/2]+1}^n M\big(Z_i-\cdot,\vec{h}\big),\;\vec{h}\in
\cH_{[n/2]}^d \bigg\}.
\end{equation}
Our  goal is to select an estimator from this collection and to study its properties. Following our general receipt we suggest  to proceed as follows.
\begin{itemize}
\item
Consider first the estimation of $A_g(f)=(1-\alpha)f+\alpha [g\star f]$, where, remind "$\star$" denotes
the convolution operator. Note that $A_g(f)$ is the density of $Z_1$ and, therefore, can be  estimated from the observation $X_1^{(n)}$
using standard kernel estimator. Introduce the collection
\begin{equation}
\label{eq2:def-K-est}
\mathbf{A}\big(\cH_{[n/2]}^d\big)=\bigg\{\widehat{A}_{\vec{h}}(\cdot)=[n/2]^{-1}\sum_{i=1}^{[n/2]}K_{\vec{h}}(Z_i-\cdot),\;\vec{h}\in\cH_{[n/2]}^d \bigg\}.
\end{equation}

\item To each estimator from this collection associate its $\bL_p$-risk ($\ell(\cdot)=\|\cdot\|_p$)
$$
\cR_{A}\big[\widehat{A}_{\vec{h}}, f\big]=\Big(\bE^{(n)}_f\Big[\big\|\widehat{A}_{\vec{h}}-A_g(f)\big\|^p_p\Big]\Big)^{1/p},\quad f\in\mP.
$$
\item Select an estimator from the collection $\mathbf{A}\big(\cH_{[n/2]}^d\big)$ in accordance with  $(\Psi_n,\ell)$-selection rule 
based on the collection of auxiliary estimators (\ref{eq:def-K-star-K}), where  $T^{(m)}=(X_1,\ldots, X_{[n/2]}), m=[n/2]$ and
    $\Psi_n=\Delta_{[n/2]}$, $\ell(\cdot)=\|\cdot\|_p$. Thus, the selection rule is
\begin{eqnarray}
\label{eq1:selection-ruleGL}
\widehat{R}_n\big(\vec{h}\big)&:=&\sup_{\vec{\eta}\in \cH^d_{[n/2]}}\Big[\big\|\widehat{A}_{\vec{h},\vec{\eta}}-\widehat{A}_{\vec{\eta}}\big)\big\|_p-
 2\Delta_{[n/2]}(\vec{\eta})\Big]_+;
\\
\label{eq2:selection-ruleGL}
\vec{\mathbf{h}}_n&:=&\inf_{\vec{h}\in \cH^d_{[n/2]}}\big[\widehat{R}_n(\vec{h})+2\Delta_{[n/2]}(\vec{h})\big].
\end{eqnarray}

\item

Choose the estimator $\widehat{B}_{\vec{\bh}_n}(\cdot)$.
\end{itemize}
The theoretical properties of the estimator $\widehat{B}_{\vec{\bh}_n}(\cdot)$ will be then deduced with help of Proposition 1, Theorem \ref{th3}(2) and Theorem \ref{th4}.

\subsection{Some bounds for the quantity  $\cE_n(\cdot,\cdot)$.}
 \label{sec:subsec-verification-B}

Introduce the following notations. Set for any $x\in\bR^d$
$$
\chi^{(n)}_{\vec{h}}(x)=\widehat{B}_{\vec{h}}-\Upsilon_{\vec{h}}(f,x),\quad \Upsilon_{\vec{h}}(f,x)=\bE^{(n)}_{2,f}\Big\{\widehat{B}_{\vec{h}}(x)\Big\}=\int_{\bR^d}M\big(z-x,\vec{h}\big)A_g(f,z)\rd z
$$
Put for brevity $\mmp=A_g(f)$. Applying Bahr-Esseen and Rosenthal inequalities, \cite{bahr},\cite{rosen}, to the sum of i.i.d. random variables $\chi^{(n)}_{\vec{h}}(x)$, we have for any $x\in\bR^d$
\begin{gather*}
\label{eq:24}
\mathbf{c^{-1}_4}\bE^{(n)}_{2,f}\Big|\chi^{(n)}_{\vec{h}}(x)\Big|^{p}\leq n^{p-1}\int_{\bR^d}\big|M\big(z-x,\vec{h}\big)\big|^p\mmp(z)\rd z,\quad p\leq 2;
\\
\mathbf{c^{-1}_4}\bE^{(n)}_{2,f}\Big|\chi^{(n)}_{\vec{h}}(x)\Big|^{p}\leq\Big(n^{-1}\int_{\bR^d}M^2\big(z-x,\vec{h}\big)\mmp(z)\rd z\Big)^{p/2}
+n^{p-1}\int_{\bR^d}\big|M\big(z-x,\vec{h}\big)\big|^p\mmp(z)\rd z,\; p>2.
\end{gather*}
Noting that in view of the Fubini theorem
$
\bE^{(n)}_{2,f}\big\|\chi^{(n)}_{\vec{h}}\big\|_p^{p}=\big\|\bE^{(n)}_{2,f}\big|\chi^{(n)}_{\vec{h}}(\cdot)\big|^{p}\big\|_1
$
one has
\begin{eqnarray*}
\label{eq:25}
\Big(\bE^{(n)}_{2,f}\big\|\chi^{(n)}_{\vec{h}}\big\|_p^{p}\Big)^{1/p}&\leq& \mathbf{c_4}n^{1/p-1}\big\|M\big(\cdot,\vec{h}\big)\big\|_p,\quad p\leq 2;
\\
\Big(\bE^{(n)}_{2,f}\big\|\chi^{(n)}_{\vec{h}}\big\|_p^{p}\Big)^{1/p}&\leq&\mathbf{c_4}n^{-\frac{1}{2}}\bigg\{\int_{\bR^d}\Big(\int_{\bR^d}M^2\big(z-x,\vec{h}\big)
\mmp(z)\rd z\Big)^{p/2}\rd x\bigg\}^{1/p}
+\mathbf{c_4}n^{\frac{1-p}{p}}\big\|M\big(\cdot,\vec{h}\big)\big\|_p
\\*[2mm]
&\leq& \mathbf{c_4}\Big(\|\mmp\|_\infty^{1/2-1/p}n^{-\frac{1}{2}}\big\|M\big(\cdot,\vec{h}\big)\big\|_2+n^{\frac{1-p}{p}}\big\|M\big(\cdot,\vec{h}\big)\big\|_p\Big),
\quad p>2
\end{eqnarray*}
To get the last estimate we have used again the Young inequality.
Thus, for any $\vec{h}\in\cH^d$ we have the following bound for the quantity $\cE(\cdot,\cdot)$ defined in (\ref{eq:def-calE})
\begin{eqnarray}
\label{eq:26}
\cE\big(f,\vec{h}\big):=\Big(\bE^{(n)}_{2,f}\big\|\chi^{(n)}_{\vec{h}}\big\|_p^{p}\Big)^{1/p}\leq \mathbf{c_5}\Big(n^{-1/2}\big\|M\big(\cdot,\vec{h}\big)\big\|_2+
n^{1/p-1}\big\|M\big(\cdot,\vec{h}\big)\big\|_p\Big), \quad\forall f\in\mP_p.
\end{eqnarray}
We would like to stress that if $\alpha=1$ and $\|g\|_\infty<\infty$ then $\mP_p=\mP$ for all $p\geq 1.$


 Let $\check{M}\big(t,\vec{h}\big), t\in\bR^d,$ denote the Fourier transform of $M\big(\cdot,\vec{h}\big)$.
Then, we obtain in view of the definition of $M\big(\cdot,\vec{h}\big)$
\begin{eqnarray}
\label{eq:fourier-trans-of-kernel-M}
\check{M}\big(t,\vec{h}\big)=\check{K}\big(t\vec{h}\big)\big[(1-\alpha)+\alpha\check{g}(-t)\big]^{-1},\;\; t\in\bR^d.&&
\end{eqnarray}
The conditions imposed on $K$ guarantee that $\check{M}\big(\cdot,\vec{h}\big)\in\bL_1\big(\bR^d\big)\cap\bL_2\big(\bR^d\big)$ for any $\vec{h}\in\cH^d$  and, hence,
$$
\big\|M\big(\cdot,\vec{h}\big)\big\|_\infty\leq (2\pi)^{-d}\big\|\check{M}\big(\cdot,\vec{h}\big)\big\|_1,\qquad
\big\|M\big(\cdot,\vec{h}\big)\big\|_2= (2\pi)^{-d}\big\|\check{M}\big(\cdot,\vec{h}\big)\big\|_2.
$$
Set  $\vec{\bmu}(\alpha)=\vec{\mu}$, $\alpha=1$, $\vec{\bmu}(\alpha)=(0,\ldots,0)$, $\alpha\in [0,1)$. We have in view of Assumption \ref{ass1:ass-on-noise-upper-bound}  for any $\vec{h}\in\cH^d$
\begin{eqnarray}
\label{eq:bound-for-supnorm-and-L_2-norm}
\big\|M\big(\cdot,\vec{h}\big)\big\|_\infty\leq \mathbf{c_6}\prod_{j=1}^dh_j^{-1-\blg_j(\alpha)},\qquad
\big\|M\big(\cdot,\vec{h}\big)\big\|_2\leq \mathbf{c_6}\prod_{j=1}^dh_j^{-\frac{1}{2}-\blg_j(\alpha)},
&&
\end{eqnarray}
Additionally we deduce from (\ref{eq:bound-for-supnorm-and-L_2-norm}) for any $p>2$
\begin{eqnarray}
\label{eq:bound-for-Lp-norm-p>2}
\big\|M\big(\cdot,\vec{h}\big)\big\|_p\leq \big\|M\big(\cdot,\vec{h}\big)\big\|^{1-\frac{2}{p}}_\infty\big\|M\big(\cdot,\vec{h}\big)\big\|_2^{\frac{2}{p}}
\leq \mathbf{c_6}\prod_{j=1}^dh_j^{-1+1/p-\blg_j(\alpha)},\quad \forall\vec{h}\in\cH^d.&&
\end{eqnarray}
If $p<2$ we will study only the case $\alpha<1/2$. Then we have in view of the definition of $M\big(\cdot,\vec{h}\big)$, applying the Young and triangle inequalities
$$
\big\|K_{\vec{h}}\big\|_p\geq (1-\alpha)\big\|M\big(\cdot,\vec{h}\big)\big\|_p-\alpha\big\|M\big(\cdot,\vec{h}\big)\big\|_p=(1-2\alpha)\big\|M\big(\cdot,\vec{h}\big)\big\|_p.
$$
It yields, for any $p<2$ and $\alpha<1/2$
\begin{eqnarray}
\label{eq:bound-for-Lp-norm-p<2}
\big\|M\big(\cdot,\vec{h}\big)\big\|_p\leq (1-2\alpha)^{-1}\|K\|_p\prod_{j=1}^dh_j^{-1+1/p},\quad \forall\vec{h}\in\cH^d.&&
\end{eqnarray}
Thus, we obtain from (\ref{eq:26}), (\ref{eq:bound-for-supnorm-and-L_2-norm}), (\ref{eq:bound-for-Lp-norm-p>2}) and (\ref{eq:bound-for-Lp-norm-p<2})
for any $\vec{h}\in\cH_{[n/2]}^d$
\begin{eqnarray}
\label{eq:27}
\sup_{f\in\mP_p}\cE\big(f,\vec{h}\big)\leq\mathbf{c_7}
\left\{\begin{array}{ll}
\; n^{-\frac{1}{2}}\prod_{j=1}^dh_j^{-\frac{1}{2}-\blg_j(\alpha)},\quad &p\geq 2,\; \alpha\in[0,1];
\\*[2mm]
n^{\frac{1}{p}-1}\prod_{j=1}^dh_j^{-1+1/p-\blg_j(\alpha)},\quad &p<2,\; \alpha\in[0,1/2).
\end{array}\right.
\end{eqnarray}
To get the bound for $p\geq 2$ we have taken into account that $nV_{\vec{h}}\geq \ln(n)$ if $\vec{h}\in\cH_{[n/2]}^d$.

In fact (\ref{eq:26}) is true for all $\alpha\in[0,1]$ but it requires to impose several additional assumption on $g$ and we do not treat the case $p<2$, $\alpha>1/2$ in the present paper. The interested reader can look the paper \cite{rebelles15}, where the corresponding norm is computed when $\alpha=1$.

\subsection{Consequences of Proposition \ref{prop1} and Theorem \ref{th3} (2). Case $\alpha\in (0,1)$.}
 \label{sec:subsec-conseq-prop1}

Let us compute the approximation error related to the estimator $\widehat{B}_{\vec{h}}$.
Note that for any $\vec{h}\in\cH^d$
\begin{eqnarray*}
\Upsilon_{\vec{h}}(f,x)&:=&\int_{\bR^d}M\big(z-x,\vec{h}\big)A_g(z)\rd z
\\
\nonumber
&=&(1-\alpha)\int_{\bR^d}M\big(z-x,\vec{h}\big)f(z)\rd z+\alpha\int_{\bR^d}M\big(z-x,\vec{h}\big)\big[f\star g\big](z)\rd z,
\\
&=&\int_{\bR^d}f(t)
\Big[(1-\alpha)M\big(t-x,\vec{h}\big)+\alpha\int_{\bR^d}M\big(u,\vec{h}\big)g(u-[t-x])\rd u\Big]\rd t
\\
&=&\int_{\bR^d}K_{\vec{h}}(t-x)f(t)\rd t.
\end{eqnarray*}
Thus, we have  for any $\vec{h}\in\cH^d$ and any $p\geq 1$
\begin{eqnarray}
\label{eq28}
\|\Upsilon_{\vec{h}}(f)-f\|_p=\|K_{\vec{h}}\star f-f\|_p.&&
\end{eqnarray}

\subsubsection{Verification of the condition (\ref{eq:hyp-cor1}) and (\ref{eq:hyp-cor2}) of Proposition \ref{prop1}.}
We will study two different cases: either $p\neq 2,\; \alpha<1/2$ or $p=2,\; \alpha\in (0,1)$.

\smallskip

Set $\Lambda_{\vec{h}}(f,\cdot)=\bE^{(n)}_{1,f}\big\{\widehat{A}_{\vec{h}}(\cdot)\big\}$, $\vec{h}\in (0,1]^d$ and note that
\begin{equation}
\label{eq:Lambda}
\Lambda_{\vec{h}}(f,\cdot)=\int_{\bR^d}K_{\vec{h}}(z-\cdot)A_g(f,z)\rd z=(1-\alpha)\big[K_{\vec{h}}\star f](\cdot)+\alpha\big[K_{\vec{h}}\star g\star f](\cdot).
\end{equation}

We have  in view of (\ref{eq:Lambda}) and the Young  inequality (recall that $g$ is a density)
\begin{eqnarray*}
\|\Lambda_{\vec{h}}(f)-A_g(f)\|_p&=&\big\|(1-\alpha)\big\{K_{\vec{h}}\star f-f\big\}+\alpha g\star\big\{K_{\vec{h}}\star f -f\big\} \big\|_p
\\
&\geq&
(1-2\alpha)\big\|K_{\vec{h}}\star f -f\big\|_p.
\end{eqnarray*}
Thus, if $\alpha<1/2$ one has for any $p\geq 1$ in view of (\ref{eq28})
$$
\|\Upsilon_{\vec{h}}(f)-f\|_p\leq (1-2\alpha)^{-1}\|\Lambda_{\vec{h}}(f)-A_g(f)\|_p
$$
and, therefore, the condition (\ref{eq:hyp-cor1}) of Proposition \ref{prop1} is fulfilled with $C_\ell=(1-2\alpha)^{-1}$.

If $p=2$ we have in view of (\ref{eq:Lambda}), (\ref{eq28}), Assumption \ref{ass1:ass-on-noise-upper-bound} and the Parseval identity
\begin{eqnarray*}
\|\Lambda_{\vec{h}}(f)-A_g(f)\|_2&=&\big\|(1-\alpha)\big[\check{K}(\vec{h}\cdot)-1\big]\check{f}(\cdot)+
\alpha\big[\check{K}(\vec{h}\cdot)-1\big]\check{g}(\cdot)\check{f}(\cdot)\big\|_2
\\*[1mm]
&=&
\big\|\big\{(1+\alpha)+\alpha\check{g}(\cdot)\big\}\big[\check{K}(\vec{h}\cdot)-1\big]\check{f}(\cdot)\big\|_2\geq \varpi
\big\|\big[\check{K}(\vec{h}\cdot)-1\big]\check{f}(\cdot)\big\|_2
\\*[1mm]
&=&\varpi\|\Upsilon_{\vec{h}}(f)-f\|_2.
\end{eqnarray*}
We conclude that the condition (\ref{eq:hyp-cor1}) of Proposition \ref{prop1} is fulfilled with $C_\ell=\varpi^{-1}$ if $p=2$.

Let us now check (\ref{eq:hyp-cor2}) of Proposition \ref{prop1}.
Note that the definition of $\Delta_{[n/2]}$ implies
\begin{eqnarray}
\label{eq29}
\Delta_{[n/2]}\big(\vec{h}\big)\geq\mathbf{c_9}
\left\{\begin{array}{ll}
\; n^{-\frac{1}{2}}\prod_{j=1}^dh_j^{-\frac{1}{2}},\quad &p\geq 2;
\\*[2mm]
n^{\frac{1}{p}-1}\prod_{j=1}^dh_j^{-1+1/p},\quad & p<2.
\end{array}\right.
\end{eqnarray}
and therefore, (\ref{eq:hyp-cor2}) is fulfilled with $C_{\cE}=\mathbf{c_9}/\mathbf{c_7}$ in view of (\ref{eq:27}), since recall  $\vec{\blg}(\alpha)=0, \alpha\neq 1$.

\subsubsection{Main results. } We deduce from Proposition \ref{prop1} and Lemma \ref{lem:extract-from-GL11}  the following assertion.
\begin{theorem}
\label{th:oracle-decon-alpha-neq-1}
Let either $p\neq 2,\; \alpha<1/2$ or $p=2,\; \alpha\in (0,1)$. Then for all $n\geq 2$
\begin{eqnarray*}
\cR_{B}\big[\widehat{B}_{\vec{\mathbf{h}}_n}, f\big]\leq \mathbf{c_{10}}\inf_{\vec{h}\in H_n}\left\{\big\|K_{\vec{h}}\star f- f\big\|_p+\big(nV_{\vec{h}}\big)^{\frac{1}{p\wedge2}-1}\right\}+\e_{[n/2]},\quad\forall f\in\mP_p.
\end{eqnarray*}
\end{theorem}
We also have the following adaptive result.
\begin{theorem}
\label{th:adaptive-decon-alpha-neq-1}
Let either $p\neq 2,\; \alpha<1/2$ or $p=2,\; \alpha\in (0,1)$. Then
the estimator $\widehat{B}_{\vec{\mathbf{h}}_n}$ is optimally-adaptive over the scale of anisotropic
 classes $\big\{\bN_{p,d}\big(\vec{\beta},L\big)\cap\mP_p,\; \vec{\beta}\in (0,\ms]^d, L>0 \Big\}$.
\end{theorem}
The fact that $\widehat{B}_{\vec{\mathbf{h}}_n}$ is $S_n$-adaptive with
$S_n=\Big\{s_{[n/2]}\big(\bN_{p,d}\big(\vec{\beta},L\big)\cap\mP_p\big),\; \vec{\beta}\in (0,\ms]^d, L>0 \big\}$
follows from (\ref{eq:23-adaptive-new}) with $\mmp=f$ and $m=[n/2]$, the  assertion (\ref{eq3:lem1}) of Lemma \ref{lem:extract-from-GL11} and the second assertion of Theorem \ref{th3}.

The  lower bound for minimax risk showing that $S_n$ is the family of minimax rates in this problem can be found in
\cite{LW15}, Theorem 1.
\begin{remark}
The results presented  in Theorems \ref{th:oracle-decon-alpha-neq-1} and \ref{th:adaptive-decon-alpha-neq-1} are new.
\end{remark}

\subsection{Consequences of  Theorem \ref{th4}. Case $\alpha=1, p=2,\; d=1$.}
 \label{sec:subsec-conseq-th4-alpha=1}
In this section we will study the adaptive estimation over the collection of univariate Sobolev classes.
\begin{definition}
\label{def:anisotrop-sobolev}
Let  $\beta_1\in (0,\infty)^d$ and $L>0$ be given. We say that $Q\in\bL_1(\bR^1)$ belongs to  Sobolev class $\mathbb{W}\big(\beta_1,L\big)$ if
$$
\int_{\bR^1}\big(1+t^2\big)^{\beta_1}\big|\check{Q}(t)\big|^2\rd t\leq L^2.
$$

\end{definition}
The adaptive estimation over the collection of {\it anisotropic} Sobolev classes when  $\alpha=1$ was studied in \cite{comte}.

\paragraph{Auxiliary inequalities.} The  results below follow from the Parseval identity,  the properties of the kernel $K$ and Assumption \ref{ass1:ass-on-noise-upper-bound}. First note, that for any $f\in\mathbb{W}\big(\beta_1,L\big)$ and any $h_1\in(0,1)$
\begin{eqnarray*}
\|K_{h_1}\star A_g(f)-A_g(f)\|_2\leq G_2Lh_1^{\beta_1+\mu_1},
\end{eqnarray*}
since obviously $A_g(f)\in \mathbb{W}\big(\beta_1+\mu_1,G_2L\big)$. 

\smallskip

It  yields in view of (\ref{eq:20}) and (\ref{eq:21}) with $\mmp=A_g(f)$ and $p=2$,

\begin{eqnarray}
\label{eq30}
s_n\big(\bF_\ma\big)&=&\sup_{f\in\bF_\ma}\inf_{\mh\in\mH_n}\big\{\cB^{(n)}_A(f,\mh)+\psi_n(f,\mh)\big\}
\nonumber\\
&=:&
\sup_{f\in\mathbb{W}(\beta_1,L)}\inf_{h_1\in H^1_n}\left\{(1+2\|K\|_1)\|K_{h_1}\star A_g(f)-A_g(f)\|_2+\mathbf{c_{11}}\big(nh_1\big)^{-\frac{1}{2}}\right\}
\nonumber\\
&\leq&\mathbf{c_{12}}L^{\frac{1}{2\beta_1+2\mu_1+1}}n^{-\frac{\beta_1+\mu_1}{2\beta_1+2\mu_1+1}}
=:\mathbf{c_{12}}
s_n\big(\mathbb{W}\big(\beta_1,L\big)\big).
\end{eqnarray}
We also have in view of Assumption \ref{ass1:ass-on-noise-upper-bound} for any $f\in\mathbb{W}\big(\beta_1,L\big)$, $y>0$ and $h_1\in(0,1]$
\begin{eqnarray*}
\|K_{h_1}\star f-f\|^2_2&\leq&\int_{-y}^y |\check{K}(h_1t)-1|^2|\check{f}(t)|^2\rd t +(1+\|K\|_1)^2\int_{\bR} \mathrm{1}_{\bR\setminus[-y,y]}(t)|\check{f}(t)|^2\rd t
\\
&\leq&\mathbf{c_{13}}\bigg[(1+y^{2\mu_1})\int_{-y}^y |g(t)|^2|\check{K}(h_1t)-1|^2|\check{f}(t)|^2\rd t +L^2y^{-2\beta_1}\bigg]
\\
&\leq&\mathbf{c_{13}}\bigg[(1+y^{2\mu_1})\|K_{h_1}\star A_g(f)-A_g(f)\|^2_2 +L^2y^{-2\beta_1}\bigg].
\end{eqnarray*}
Minimizing the r.h.s. of the latter inequality in $y$ we get for any $f\in\mathbb{W}\big(\beta_1,L\big)$ and $h_1\in(0,1]$
\begin{equation}
\label{eq31}
\|K_{h_1}\star f-f\|_2\leq\mathbf{c_{13}}\bigg[\|K_{h_1}\star A_g(f)-A_g(f)\|_2+L^{\frac{\mu_1}{\mu_1+\beta_1}}\|K_{h_1}\star A_g(f)-A_g(f)\|_2^{\frac{\beta_1}{\beta_1+\mu_1}}\bigg].
\end{equation}

\paragraph{Application of Theorem \ref{th4}.}

First, we note that all assumptions of Theorem \ref{th4} are fulfilled. Indeed, in Section \ref{sec:subsec-verification-A} we have already checked  the hypotheses ${\bf A^{\text{permute}}}$,
${\bf A^{\text{upper}}}$ and ${\boldsymbol \cA^{\text{upper}}}$. Moreover, defined in (\ref{eq:def-psi_h(F_ma)}) quantity in our case is given by
$$
\psi_n\big(\mathbb{W}\big(\beta_1,L\big),h_1\big)=\mathbf{c_{14}}\big(nh_1\big)^{-\frac{1}{2}}.
$$
It yields together with (\ref{eq30}) that
$$
\cV_n(\ma)=\Big\{h_1\in\cH^1_n:\; h_1\geq \mathbf{c_{15}}L^{-\frac{2}{2\beta_1+2\mu_1+1}}n^{-\frac{1}{2\beta_1+2\mu_1+1}}\Big\},\quad \ma=(\beta_1,L),
$$
and we have that in view of (\ref{eq:27})
\begin{equation}
\label{eq32}
\sup_{f\in\bF_\ma}\sup_{\mh\in\cV_n(\ma)}\cE_n\big(f,h_1\big)\leq \mathbf{c_{16}}L^{\frac{1+2\mu_1}{2\beta_1+2\mu_1+1}}n^{-\frac{\beta_1}{2\beta_1+2\mu_1+1}}.
\end{equation}
Since $\cU_n(\ma,f)=\big\{h_1\in\cH^1_n:\; \|K_{h_1}\star A_g(f)-A_g(f)\|_2< 8s_n\big(\mathbb{W}\big(\beta_1,L\big)\big)\big\}$
we obtain from  (\ref{eq31})
\begin{equation*}
\sup_{f\in\bF_\ma}\sup_{\mh\in\cU_n(\ma,f)}\|K_{h_1}\star f-f\|_2\leq
\mathbf{c_{17}}L^{\frac{\mu_1}{\mu_1+\beta_1}}\Big[s_n\big(\mathbb{W}\big(\beta_1,L\big)\big)\Big]^{\frac{\beta_1}{\beta_1+\mu_1}}
=\mathbf{c_{17}}L^{\frac{1+2\mu_1}{2\beta_1+2\mu_1+1}}n^{-\frac{\beta_1}{2\beta_1+2\mu_1+1}}.
\end{equation*}
It yields  together with  (\ref{eq32})
\begin{equation}
\label{eq33}
\varphi_n\big(\mathbb{W}\big(\beta_1,L\big)\big)=\mathbf{c_{18}}L^{\frac{1+2\mu_1}{2\beta_1+2\mu_1+1}}n^{-\frac{\beta_1}{2\beta_1+2\mu_1+1}}.
\end{equation}
where $\mathbf{c_{18}}$ independent of $L$.  Putting  $\mu^{*}=\max_{j=1,\ldots,d}\mu_j$ we note that
$$
\delta_n\geq \mathbf{c_{19}}n^{-\frac{1}{2}},\quad \sup_{f\in\mP}\sup_{h_1\in\cH^1_n}\cE_n(f,h_1)\leq \mathbf{c_{20}}n^{\mu^*},\;
\sup_{f\in\mathbb{W}\big(\beta_1,L\big)}\sup_{h_1\in (0,1]}\|K_{h_1}\star f-f\|_2\leq 2\|K\|_1L.
$$
Hence, in view of the assertion (\ref{eq3:lem1}) of Lemma \ref{lem:extract-from-GL11} we can state that
\begin{equation}
\label{eq333}
\limsup_{n\to\infty}n^{a}\kappa_n\big(\mathbb{W}\big(\beta_1,L\big)\big)=0, \quad \forall a>0.
\end{equation}
Thus, we deduce from (\ref{eq31}), (\ref{eq33}), (\ref{eq333})  and Theorem \ref{th4} the following result.
\begin{assertion}
\label{asser1}
Let $\mathbf{h}_n$ comes from (\ref{eq1:selection-ruleGL})--(\ref{eq2:selection-ruleGL}). Then for any $\mu_1>0,$ $\beta_1>0$ and $L>0$
\begin{eqnarray*}
\limsup_{n\to\infty}\varphi^{-1}_n\big(\mathbb{W}\big(\beta_1,L\big)\big)\sup_{f\in\mathbb{W}\big(\beta_1,L\big)}\Big\{\bE^{(n)}_{f}\big\|
\widehat{B}_{\vec{\mathbf{h}}_n}- f\big\|^2_2\Big\}^{\frac{1}{2}}<\mathbf{c_{21}},
\end{eqnarray*}
where $\mathbf{c_{21}}$ is independent of $L$.

The estimator $\widehat{B}_{\vec{\mathbf{h}}_n}$ is optimally-adaptive over the scale $\big\{\mathbb{W}\big(\beta_1,L\big), \beta_1>0, L>0\big\}$.
\end{assertion}
Assertion \ref{asser1} is the particular case of the results obtained in \cite{comte}, which were established by the use of completely different selection scheme.

\section{Simultaneous adaptive estimation of partial derivatives in the density model.}
\label{sec:derivatives}

Let as previously $X_i, i=1,\ldots,n,$ are {\it i.i.d.} $d$-dimensional  random vectors with unknown density $f$.
This is the particular case of the model considered above corresponding  to $\alpha=0$.

For any multi-index $\mm=(m_1,\ldots,m_d)\in\bN^d$ let
$
f^{(\mm)}=\frac{\partial^{m_1+\cdots +m_d}f}{\partial x_1^{m_1}\cdots \partial x_d^{m_d}}
$
be the partial derivative of $f$ of order $\mm$. The goal is to estimate
$B(f)=f^{(\mm)}$ under  $\bL_p$-loss, that is 
$$
\cR_{B}\big[\widetilde{B}_n, f\big]=\Big(\bE^{(n)}_f\Big[\big\|\widetilde{B}_n-f^{(\mm)}\big\|^p_p\Big]\Big)^{1/p},\quad f\in\mP.
$$
\begin{remark}
We will keep all previous notations with only one change related to the following simple observation. All oracle results established in the present paper
remain valid if one replace in all formulas the set $\mH_n$ (abstract model) or $\cH^d_m$ (density model) by its arbitrary subset.
In particular we will use below the selection rule (\ref{eq1:selection-ruleGL})--(\ref{eq2:selection-ruleGL})
from the collection $\mathbf{A}(\cdot)$, which is parameterized by  either $\cH_{[n/2]}^d$  as it was before or by its subset
$\cH^{\text{isotrop}}_n:=\big\{\vec{h}\in\cH_{[n/2]}^d : h_1=\cdots=h_d\big\}$.
The latter set will be used when the adaptation over the scale of isotropic,
 ($\beta_1=\cdots=\beta_d$)
Nikol'skii classes is studied.
\end{remark}
Thus, let $H_n$ denote either $\cH_{[n/2]}^d$ or $\cH^{\text{isotrop}}_n$ and let the collection
$\mathbf{A}\big(H_n\big)$ be given by (\ref{eq2:def-K-est}) (where $\cH_{[n/2]}^d$ is replaced by $H_n$) and  $Z_i=X_i, i=1,\ldots,[n/2]$.
Introduce the family of estimators

$$
\mathbf{B}\big(H_n\big)=\bigg\{
\widehat{B}_{\vec{h}}(\cdot)=(n-[n/2])^{-1}\sum_{i=[n/2]+1}^n \prod_{j=1}^d(-1)^{m_j}h_j^{-m_j}
K^{(\mm)}_{\vec{h}}\big(X_i-\cdot\big),\;\;\vec{h}\in H_n\bigg\},
$$
where $K^{(\mm)}$ denotes the partial derivative of $K$ of the order $\mm$.

It is important to note that the estimator $\widehat{B}_{\vec{h}_n}$ would coincide with $\widehat{A}^{(\mm)}_{\vec{h}_n}$ if we would construct it
using observations $X_i, i=1,\ldots,[n/2]$ instead of $X_i, i=[n/2]+1,\ldots,n$.
The interest to this collection is dictated by the following minimax result. For
 any $\vec{\beta}\in (0,\infty)^d$ and $\mm\in\bN^d$ put
 $$
 \frac{1}{\beta}=\sum_{j=1}^d\frac{1}{\beta_j}, \qquad
 \frac{1}{\omega}=\sum_{j=1}^d\frac{m_j}{\beta_j}.
 $$
\begin{proposition}
\label{prop2}
For any $p\geq 1, L>0$, any $\vec{\beta}\in (0,\infty)^d$ and $\mm\in\bN^d$ provided $\omega >1$ one can find a kernel $K$ and 
$\vec{h}_n\in (0,1]^d$  such that the estimator $\widehat{A}^{(\mm)}_{\vec{h}_n}$ is rate-minimax for $f^{(\mm)}$ on
$\bN_{p,d}\big(\vec{\beta},L\big)\cap\mP_p$.

The minimax rate of convergence is given by
$$
\varphi^{(\mm)}_n\big(\bN_{p,d}\big(\vec{\beta},L\big)\cap\mP_p\big)=L^{\frac{(1/\beta)(1-1/p\wedge 2)+1/\omega}{1+(1/\beta)(1-1/p\wedge 2)}}
n^{-\frac{(1-1/\omega)(1-1/p\wedge 2)}{1+(1/\beta)(1-1/p\wedge 2)}}.
$$
The conditions $\omega>1$ and $p>1$  are necessary for the existence of uniformly consistent on $\bN_{p,d}\big(\vec{\beta},L\big)$ estimators.

\end{proposition}

The assertions of the proposition seems to be new.

Now, let us formulate the main result of this section.
As previously let $\ms\in\bN^*$ be an arbitrary but a priory chosen number. Let the kernel $K$, involved to the description of estimators from the collections
$\mathbf{A}\big(H_n\big)$ and $\mathbf{B}\big(H_n\big)$, be constructed in accordance with (\ref{eq:w-function}).

Set at last
$
\Pi_\ms=\big\{(\vec{\beta},\mm)\in (0,\ms]^d\times\bN^d:\; \omega>1\big\}.
$

\begin{theorem}
\label{th7}
Let $\vec{\mathbf{h}}_n$ be  the issue of (\ref{eq1:selection-ruleGL})--(\ref{eq2:selection-ruleGL}) with $H_n=\cH_{[n/2]}^d$ if $m_1=\cdots=m_d$
or with $H_n=\cH^{\text{isotrop}}_n$ if $\beta_1=\cdots=\beta_d$.
Then for any $p\geq 1$, $L>0,$ and $(\vec{\beta},\mm)\in\Pi_\ms$
\begin{eqnarray*}
\limsup_{n\to\infty}\Big[\varphi^{(\mm)}_n\big(\bN_{p,d}\big(\vec{\beta},L\big)\cap\mP_p\big)\Big]^{-1}
\sup_{f\in\bN_{p,d}\big(\vec{\beta},L\big)\cap\mP_p}\Big\{\bE^{(n)}_{f}\big\|
\widehat{B}_{\vec{\mathbf{h}}_n}- f\big\|^p_p\Big\}^{\frac{1}{p}}<\mathbf{c_{21}},
\end{eqnarray*}
where  $\mathbf{c_{21}}$ is independent of $L$.
\end{theorem}

Theorem \ref{th7} together with second statement  of Proposition \ref{prop2} allows us to assert that  the estimator
$\widehat{B}_{\vec{\mathbf{h}}_n}$ is optimally-adaptive over the scale of either anisotropic Nikol'skii classes (if $m_1=\cdots=m_d$) or isotropic ones
(without any restriction on the order of the considered partial derivative). Note that the problem is completely solved in the dimension 1.
Also we conclude that the differentiation of an optimally-adaptive estimator leads to the optimally-adaptive
estimator of the corresponding
partial derivative.

\subsection{Proof of Proposition \ref{prop2}.}

\begin{lemma}
\label{lem:-embedding-for-derivatives}
For any $p\geq 1$,  any $\vec{\beta}$ and  $\mm$, provided $\omega>1$,
and any $f\in\bN_{p,d}\big(\vec{\beta},L\big)$ one has
\begin{equation}
\label{eq:bound-for-derivative}
\|f^{(\mm)}\|_p\leq \mathbf{c_{23}}L^{1/\omega}\|f\|_p^{1-1/\omega},
\end{equation}
where $\mathbf{c_{23}}$ is independent of $L$.

\end{lemma}
The author persuades  that this {\it Kolmogorov-type} inequality should be known but he was unable to find the exact reference.


\paragraph{Proof of Lemma \ref{lem:-embedding-for-derivatives}.} The following inequality
 was proved in \cite{GL13}.
For any $p\geq 1$ and  $\ms\in\bN^d$ there exists $\mathbf{c_{22}}$ is independent of $L$ such that
\begin{equation}
\label{eq:bound-for bias}
\sup_{F\in\bN_{p,d}\big(\vec{\beta},L\big)}\|K_{\vec{h}}\star F-F\|_p\leq \mathbf{c_{22}}L\sum_{j=1}^dh_j^{\beta_j},
\quad\forall \vec{h}\in (0,1]^d,\;\forall\vec{\beta}\in (0,\ms]^d.
\end{equation}
The following statement can be found in \cite{Nikolski}, Chapter 5, Theorem 5.6.3.

If $\omega>1$ and $f\in\bN_{p,d}\big(\vec{\beta},L\big)$ then
$f^{(\mm)}$ exists and
\begin{equation}
\label{eq:nikolski-inclusion-for-derivative}
f^{(\mm)}\in \bN_{p,d}\big(\vec{\gamma},\mathbf{c}L\big),\quad \gamma_j=\beta_j(1-1/\omega),\; j=1,\ldots,d,
\end{equation}
and $\mathbf{c}>0$ is independent of $L$.

Since $\vec{\beta}$ is fixed we can always choose $\ms\in\bN^d$ in order to have $\vec{\beta}\in (0,\ms]^d$ and, then, (\ref{eq:bound-for bias})
will be fulfilled. We obviously have by the triangle inequality
$$
\|f^{(\mm)}\|_p\leq \Big\|\big(K_{\vec{h}}\big)^{(\mm)}\star f- f^{(\mm)}\Big\|_p+\Big\|\big(K_{\vec{h}}\big)^{(\mm)}\star f\Big\|_p.
$$
It is easy to see that
$
\big(K_{\vec{h}}\big)^{(\mm)}=\prod_{j=1}^d(-1)^{m_j}h_j^{-m_j}
K^{(\mm)}_{\vec{h}_n}
$
for any $\vec{h}\in (0,1]^d$ and we have first by integrating by parts
\begin{equation}
\label{eq500}
\big(K_{\vec{h}}\big)^{(\mm)}\star f- f^{(\mm)}=K_{\vec{h}}\star f^{(\mm)}- f^{(\mm)}.
\end{equation}
Hence, applying (\ref{eq:bound-for bias}) with $F=f^{(\mm)}$ and $\vec{\beta}$ replaced by $\vec{\gamma}$
we obtain in view of (\ref{eq:nikolski-inclusion-for-derivative})
\begin{equation}
\label{eq50}
\Big\|\big(K_{\vec{h}}\big)^{(\mm)}\star f- f^{(\mm)}\Big\|_p\leq \mathbf{c_{24}}L\sum_{j=1}^dh_j^{\gamma_j},
\quad\forall \vec{h}\in (0,1]^d.
\end{equation}
Moreover, we obtain applying the Young inequality
$$
\Big\|\big(K_{\vec{h}}\big)^{(\mm)}\star f\Big\|_p\leq \Big\|\big(K_{\vec{h}}\big)^{(\mm)}\Big\|_1\|f\|_p=\|f\|_p\big\|K^{(\mm)}\big\|_1\prod_{j=1}^dh_j^{-m_j},
$$
that yields together with (\ref{eq50}) for any $f\in\bN_{p,d}\big(\vec{\beta},L\big)$
$$
\|f^{(\mm)}\|_p\leq \mathbf{c_{25}}\bigg\{L\sum_{j=1}^dh_j^{\beta_j(1-1/\omega)}+\|f\|_p\prod_{j=1}^dh_j^{-m_j}\bigg\},
 \quad\forall \vec{h}\in (0,1]^d.
$$
Choosing $h_j=(||f||_p/L)^{1/\beta_j}$ (that is possible since $||f||_p\leq L$ in view of the definition of the Nikol'skii class)
we come to the assertion of the lemma.
\epr

\noindent{\bf Proof of the proposition.} Since
$$
\widehat{A}^{(\mm)}_{\vec{h}_n}(\cdot)=[n/2]^{-1}\sum_{i=1}^{[n/2]}\prod_{j=1}^d(-1)^{m_j}h_j^{-m_j}
K^{(\mm)}_{\vec{h}_n}\big(X_i-\cdot\big)
$$
we have in view of  (\ref{eq50}) for any $p\geq 1$, $f\in\bN_{p,d}\big(\vec{\beta},L\big)$ and  $\vec{h}\in (0,1]^d$
\begin{equation}
\label{eq51}
\Big\|\bE_{1,f}^{(n)}\big(\widehat{A}^{(\mm)}_{\vec{h}}\big)-f^{(\mm)}\Big\|_p=
\Big\|\big(K_{\vec{h}}\big)^{(\mm)}\star f- f^{(\mm)}\Big\|_p\leq \mathbf{c_{24}}L\sum_{j=1}^dh_j^{\beta_j(1-1/\omega)}.
\end{equation}
Next, repeating the computations led to (\ref{eq:27}) (remind that our considerations here correspond to the case $\alpha=0$) we get
for any $p\geq 1$, $f\in\mP$ and   $\vec{h}\in (0,1]^d$
\begin{equation}
\label{eq52}
\Big[\bE_{1,f}^{(n)}\Big\|\widehat{A}^{(\mm)}_{\vec{h}}-\bE_{1,f}^{(n)}\big(\widehat{A}^{(\mm)}_{\vec{h}}\big)\Big\|^p_p\Big]^{\frac{1}{p}}\leq
\mathbf{c_{26}}\prod_{j=1}^dh_j^{-m_j}\Big(n\prod_{j=1}^dh_j\Big)^{\frac{1}{p\wedge 2}-1}.
\end{equation}
We deduce from (\ref{eq51}) and (\ref{eq52}) that for any $p\geq 1$, $f\in\bN_{p,d}\big(\vec{\beta},L\big)\cap\mP_p$ and  $\vec{h}\in (0,1]^d$
\begin{equation*}
\cR_{A}\big[\widehat{A}_{\vec{h}}, f^{(\mm)}\big]:=
\Big[\bE_{1,f}^{(n)}\Big\|\widehat{A}^{(\mm)}-f^{(\mm)}\Big\|^p_p\Big]^{\frac{1}{p}}\leq
\mathbf{c_{27}}\bigg\{L\sum_{j=1}^dh_j^{\beta_j(1-1/\omega)}+\prod_{j=1}^dh_j^{-m_j}\Big(n\prod_{j=1}^dh_j\Big)^{\frac{1}{p\wedge 2}-1}\bigg\}.
\end{equation*}
Noting that the right hand side of the obtained inequality is independent of $f$ and  minimizing it with respect to $\vec{h}$ we come to the following bound.
\begin{equation*}
\sup_{f\in\bN_{p,d}\big(\vec{\beta},L\big)}\cR_{A}\big[\widehat{A}_{\vec{h}_n}, f^{(\mm)}\big]\leq
\mathbf{c_{28}}\varphi^{(\mm)}_n\big(\bN_{p,d}\big(\vec{\beta},L\big)\cap\mP_p\big),
\end{equation*}
where $\vec{h}_n$ is the minimizer of the right hand side of penultimate inequality.

Thus, we have proved that the maximal risk is upper-bounded by $\varphi^{(\mm)}_n\big(\bN_{p,d}\big(\vec{\beta},L\big)\cap\mP_p\big)$. The proof
of the corresponding
lower bound estimate (including the third assertion of the proposition) follows immediately from general lower bound construction established  in \cite{GL13},
Theorem 3 (tail and dense zones) and  it can be omitted.
\epr

\subsection{Proof of Theorem \ref{th7}.} The proof consists in the application of Theorem \ref{th4}. Since in Section \ref{sec:subsec-verification-A}
we have already checked  the hypothesis  ${\bf A^{\text{permute}}}$,  ${\bf A^{\text{upper}}}$ and ${\boldsymbol \cA^{\text{upper}}}$, it remains to compute
the quantities
\begin{eqnarray*}
\varphi_n(\bF_\ma)&=&\sup_{f\in\bF_\ma}\sup_{\mh\in\cV_n(\ma)}\cE_n(f,\mh)+
\sup_{f\in\bF_\ma}\sup_{\mh\in\cU_n(\ma,f)}\rho\big(\Upsilon_{\mh},B(f)\big),
\\
\kappa_n(\bF_\ma)&=&\big(5\e_n/\delta_n)\sup_{f\in\bF_\ma}\Big[
\sup_{\mh\in\mH_n}\cE_n(f,\mh)+6\sup_{\mh\in\mH_n}\rho\big(\Upsilon_{\mh},B(f)\big),
\Big],
\end{eqnarray*}
where $\cU_n(\ma,f)=\big\{\mh: \ell\big(\Lambda_{\mh},A(f)\big)< 8C_{\Psi}s_n(\bF_\ma)\big\}$ and
$
\cV_n(\ma)=\big\{\mh: \psi_n(\bF_\ma,\mh)< 4C^2_{\Psi}s_n(\bF_\ma)\big\}$.

\smallskip

In the considered problem $\bF_\ma=\bN_{p,d}\big(\vec{\beta},L\big)\cap\mP_p, \ma=\big(\vec{\beta},L\big)$, $\mh=\vec{h}\in H_n$ and
\begin{eqnarray}
\label{eq54}
\ell\big(\Lambda_{\mh},A(f)\big)=\Big\|K_{\vec{h}}\star f- f\Big\|_p,
\quad \rho\big(\Upsilon_{\mh},B(f)\big)=\Big\|\big(K_{\vec{h}}\big)^{(\mm)}\star f- f^{(\mm)}\Big\|_p
\end{eqnarray}
We also have from  (\ref{eq:21}) and (\ref{eq52})  that for any $f\in\mP_p$ and  independently of  $\vec{\beta}$ and $L$
\begin{eqnarray}
\label{eq55}
\psi_n(\bF_\ma,\mh)=\underline{\mathbf{C}}\big(nV_{\vec{h}}\big)^{\frac{1}{p\wedge2}-1},\quad \cE_n(f,\vec{h})\leq \mathbf{c_{26}}\prod_{j=1}^dh_j^{-m_j}
\Big(n\prod_{j=1}^dh_j\Big)^{\frac{1}{p\wedge 2}-1}.
\end{eqnarray}
To get the second inequality we have taken into account that the estimators $\widehat{B}_{\vec{h}}$ and
$\widehat{A}^{(\mm)}_{\vec{h}}$ have the same distribution.
Moreover 
$$
s_{[n/2]}\big(\bN_{p,d}\big(\vec{\beta},L\big)\cap\mP_p\big)=\mathbf{c_{29}}L^{\frac{1-1/(p\wedge 2)}{\beta+1-1/(p\wedge 2)}}
\;n^{-\frac{1-1/(p\wedge 2)}{1+(1/\beta)(1-1/(p\wedge 2)}}
$$ 
as it was found in (\ref{eq:23-adaptive-new}).
It remains to note that $\delta_n\geq \mathbf{c_{30}}n^{-\frac{1}{2}}$ in view of (\ref{eq:21}) and the definition of $\cH^d_{[n/2]}$ and
 moreover in view of
(\ref{eq54}), (\ref{eq55}) and (\ref{eq:nikolski-inclusion-for-derivative}) and  (\ref{eq500})
$$
\sup_{\vec{h}\in\cH^d_{[n/2]}}\sup_{f\in\mP_p}\cE_n(f,\mh)\leq n^{m^*},\quad
\sup_{f\in\bN_{p,d}\big(\vec{\beta},L\big)}\Big\|\big(K_{\vec{h}}\big)^{(\mm)}\star f- f^{(\mm)}\Big\|_p\leq \mathbf{c_{31}}L,
$$
where $m^*=\max_{j=1,\ldots,d}m_j$. To get the last bound we have also used the Young inequality.

Since the hypotheses ${\bf A^{\text{upper}}}$ and ${\boldsymbol \cA^{\text{upper}}}$ are checked with  $\e_{[n/2]}(p)$ verifying the assertion (\ref{eq3:lem1}) of
Lemma \ref{lem:extract-from-GL11} we can first assert that
\begin{eqnarray}
\label{eq56}
\limsup_{n\to\infty}n^{a}\kappa_n(\bF_\ma)=0,\quad\forall a>0.
\end{eqnarray}
Next, in view of (\ref{eq500}) we have for any $\vec{h}\in (0,1]^d$
$$
\big(K_{\vec{h}}\big)^{(\mm)}\star f- f^{(\mm)}=K_{\vec{h}}\star f^{(\mm)}- f^{(\mm)}=\big(K_{\vec{h}}\star f-f\big)^{(\mm)}
$$
and, applying Lemma \ref{lem:-embedding-for-derivatives}, we obtain for any $f\in\bN_{p,d}\big(\vec{\beta},L\big)$
\begin{eqnarray*}
\rho\big(\Upsilon_{\mh},B(f)\big)&:=&\big\|\big(K_{\vec{h}}\big)^{(\mm)}\star f- f^{(\mm)}\big\|_p\leq \mathbf{c_{23}}L^{1/\omega}
\big\|K_{\vec{h}}\star f-f\big\|_p^{1-1/\omega}
\\*[2mm]
&=:&\mathbf{c_{23}}L^{1/\omega}\ell^{1-1/\omega}\big(\Lambda_{\mh},A(f)\big).
\end{eqnarray*}
It yields for any $(\vec{\beta},\mm)\in\Pi_\ms$ and any $p\geq 1$
\begin{eqnarray}
\label{eq57}
\sup_{f\in\bF_\ma}\sup_{\mh\in\cU_n(\ma,f)}\rho\big(\Upsilon_{\mh},B(f)\big)&\leq& \mathbf{c_{32}}L^{1/\omega}
\Big[s_n\big(\bN_{p,d}\big(\vec{\beta},L\big)\cap\mP_p\big)\Big]^{1-1/\omega}
\nonumber\\
&=&
\mathbf{c_{33}}L^{\frac{(1/\beta)(1-1/p\wedge 2)+1/\omega}{1+(1/\beta)(1-1/p\wedge 2)}}
\;n^{-\frac{(1-1/\omega)(1-1/(p\wedge 2))}{1+(1/\beta)(1-1/(p\wedge 2)}}
\nonumber\\*[2mm]
&=&\mathbf{c_{33}}\varphi^{(\mm)}_n\big(\bN_{p,d}\big(\vec{\beta},L\big)\cap\mP_p\big).
\end{eqnarray}
 Put for brevity $s_n=s_n\big(\bN_{p,d}\big(\vec{\beta},L\big)\cap\mP_p\big)$, $r=p\wedge2$ and remark that
$$
\cV_n(\ma)=\Big\{\vec{h}\in H_n:\; V_{\vec{h}}^{1-\frac{1}{r}}> n^{\frac{1}{r}-1} (\mathbf{c_{33}}s_n)^{-1}
\Big\}=\Big\{\vec{h}\in H_n:\; V_{\vec{h}}> \mathbf{c_{34}}\big(Ln^{1-1/r}\big)^{-\frac{1}{\beta+1-1/r}}
\Big\}.
$$
We obtain in view of (\ref{eq55})
$$
\cE_n(f,\vec{h})\leq \mathbf{c_{26}}\left\{
\begin{array}{ll}
\;\;n^{1/r-1}V_{\vec{h}}^{1/r-1-m_1}, &\quad m_1=\cdots=m_d;
\\*[2mm]
n^{1/r-1}h_1^{d(1/r-1-d^{-1}\sum_{j=1}^dm_j)}, &\quad \vec{h}\in\cH^{\text{isotrop}}_n.
\end{array}
\right.
$$
Note that  both bounds coincide since $V_{\vec{h}}=h_1^d$ if $\vec{h}\in\cH^{\text{isotrop}}_n$ and  $m_1=d^{-1}\sum_{j=1}^dm_j$ in the first case.
Simple algebra shows that
\begin{eqnarray}
\label{eq58}
\sup_{f\in\bN_{p,d}\big(\vec{\beta},L\big)\cap\mP_p}\sup_{\vec{h}\in\cV_n(\ma)}\cE_n(f,\mh)\leq
\mathbf{c_{35}}\varphi^{(\mm)}_n\big(\bN_{p,d}\big(\vec{\beta},L\big)\cap\mP_p\big).
\end{eqnarray}
The assertion of the theorem follows now from (\ref{eq56}), (\ref{eq57}), (\ref{eq58}) and Theorem \ref{th4}.

\epr

\bibliographystyle{agsm}

\end{document}